\documentclass{article}
\pdfoutput=1

\usepackage[utf8]{inputenc}
\usepackage[a4paper, bottom = 3.5cm]{geometry}
\usepackage{amsmath}
\usepackage{amssymb}
\usepackage{bm}
\usepackage{mathtools}
\usepackage{microtype}
\usepackage{xcolor}
\usepackage[normalem]{ulem} 
\usepackage{url}
\usepackage[hidelinks]{hyperref}
\usepackage{graphicx}
\usepackage{epstopdf}
\usepackage{algorithm}%
\usepackage{algorithmicx}%
\usepackage{algpseudocode}%

\DeclareMathOperator{\spann}{span}

\algnewcommand{\LineComment}[1]{\Statex \(\triangleright\) #1}
\algrenewcommand\algorithmicrequire{\textbf{Data:}}

\title{Dynamical low-rank approximation of the Vlasov--Poisson equation with piecewise linear spatial boundary}
\author{
Andr\'e Uschmajew\thanks{Institute of Mathematics \& Centre for Advanced Analytics and Predictive Sciences, University of Augsburg, 86159 Augsburg, Germany} 
\and Andreas Zeiser\thanks{Faculty 1: School of Engineering -- Energy and Information, HTW Berlin -- University of Applied
Sciences, 12459 Berlin, Germany}
}

\date{}

\begin{document}

\maketitle

\begin{abstract}
Dynamical low-rank approximation (DLRA) for the numerical simulation of Vlasov--Poisson equations is based on separation of space and velocity variables, as proposed in several recent works. The standard approach for the time integration in the DLRA model uses a splitting of the tangent space projector for the low-rank manifold according to the separated variables. It can also be modified to allow for rank-adaptivity. A less studied aspect is the incorporation of boundary conditions in the DLRA model. In this work, a variational formulation of the projector splitting is proposed which allows to handle inflow boundary conditions on spatial domains with piecewise linear boundary. Numerical experiments demonstrate the principle feasibility of this approach.
\end{abstract}

\section{Introduction}

In this work we consider the dynamical low-rank approximation of the Vlasov--Poisson equation in $d \le 3$ spatial dimensions:
\begin{align} \label{eq:transport}
 \partial_t f + \bm v \cdot \nabla_{\bm x} f - \bm E(t, \bm x) \cdot \nabla_{\bm v} f 
	= 0 \quad \textnormal{ in } \Omega = \Omega_x \times \Omega_v, \quad t \in (0,T),
\end{align}
with a bounded domain $\Omega_x \subset \mathbb R^d$, and $\Omega_v = \mathbb R^d$. This equation models the time evolution of an electron density  $f = f(t, \bm x, \bm v)$ of a collisionless plasma as a function of space and velocity in the presence of an electrical field $\bm E$. We assume the initial condition
\[ f(0,\bm x, \bm v) = f_0(\bm x, \bm v) \]
and inflow boundary conditions of the form
\begin{equation} \label{eq:boundary_condition}
f(t,\cdot,\cdot) = g(t,\cdot,\cdot), \quad \textnormal{ on } \Gamma^- = \{(\bm x,\bm v) \in \partial \Omega_x \times  \Omega_v \,|\, \bm v \cdot \bm n_x < 0 \}.    
\end{equation}
Here $\bm n_x$ denote the outward normal vectors of the spatial domain $\Omega_x$. The electrical field $\bm E$ can either be fixed or dependent on the density $f$ via a Poisson equation:
\begin{equation} \label{eq:poisson}
    \bm E(t,\bm x) = - \nabla_{\bm x} \Phi(t,\bm x), \quad -\Delta \Phi = \rho, \quad \rho(t, \bm x) = \rho_b -\int_{ \Omega_v} f(t,\bm x, \bm v) \, \mathrm d \bm v,
\end{equation}
supplemented with appropriate boundary conditions. Here $\rho_b$ is a background charge.

Simulations of such systems are computationally demanding since the time evolution of an $2d$-dimensional, i.e.~up to six-dimensional, function has to be  calculated. Applying standard discretization schemes thus leads to an evolution equation in $\mathcal O(n^{2d})$ degrees of freedom, where $n$ is the number of grid points in one dimension and a corresponding high computational effort. To tackle the problem, methods such as  particle methods \cite{verboncoeur2005},  adaptive multiscale methods \cite{deriaz2018}, and sparse grids \cite{huang2023} have been used.

In the seminal paper~\cite{einkemmer2018a} dynamical low-rank approximation (DLRA) has been proposed for solving~\eqref{eq:transport}. DLRA is a general concept of approximating the evolution of time-dependent multivariate functions using a low-rank model, typically based on separation of variables. It originates in classic areas of mathematical physics, such as molecular dynamics~\cite{Lubich08}, and has been proposed in the works~\cite{Koch07,Koch10} as a general numerical tool for the time-integration of ODEs on fixed-rank matrix or tensor manifolds. An overview from the perspective of numerical analysis and further references can be found in~\cite{Uschmajew20}. The approach can also be made rigorous for low-rank functions in $L_2$-spaces thanks to their tensor product Hilbert space structure; see,~e.g.,~\cite{Bachmayr2021}.

For the simulation of~\eqref{eq:transport} using DLRA, a rather natural separation of space and velocity variables is applied. After discretizing the problem in a corresponding tensor product discretization space, one then seeks an approximate solution curve for equation~\eqref{eq:transport} of the form
\begin{equation}\label{eq:dlra}
f_r(t,\bm x, \bm v) = \sum_{i=1}^r\sum_{j=1}^r X_i(t,\bm x) S_{ij}(t) V_j(t,\bm v),
\end{equation}
that is, $f_r(t,\cdot,\cdot)$ is a rank-$r$ function for every time point $t$. For the numerical time-integration of the low-rank factors in this representation,~\cite{einkemmer2018a} adopted the so-called projector-splitting approach from~\cite{Lubich2014}, which is one of the work-horse algorithms for DLRA, to the case of the transport equation~\eqref{eq:transport}.

In the present work, we wish to study in more detail how to incorporate the inflow boundary condition~\eqref{eq:boundary_condition} into such a scheme. In~\cite{einkemmer2018a} this question was somewhat circumvented by considering periodic boundary conditions, which, however, is not always applicable in real world problems. Enforcing boundary conditions directly on low-rank representations like \eqref{eq:dlra} appears to be rather impractical, or at least poses some difficulties. For the nonlinear Boltzmann equation such an approach has been considered in~\cite{Hu2022}. Here, we will instead start out from a weak formulation of the transport problem~\eqref{eq:transport} including a weak formulation of the boundary condition according to \cite[Sec.~76.3]{ern_2021}. We then use the projector splitting approach for the time-stepping in this weak formulation, where the problem is iteratively reduced to subspaces belonging to variations of space or velocity variables only. We note that in~\cite{hauck2022} boundary conditions for a two-dimensional product domain (corresponding to $d=1$ in our notation) were included in this way through discretization by a discontinuous Galerkin method.

It is important to note already here that in order to still benefit from the seperation of variables, our approach assumes the boundary function $g$ to be a finite sum of tensor products, or at least be approximable by such, see~\eqref{eq:tensor_g} below. Moreover, we require the spatial domain to have a piecewise linear boundary so that the outer normal vectors are piecewise constant, see~\eqref{eq:boundary_decomposition}.

As will be derived in section~\ref{sec:dlra}, the resulting effective equations for the low-rank factors in our weak formulation of the projector splitting scheme take the form of Friedrichs' systems, that is, systems of hyperbolic equations in weak formulation, that respect the boundary conditions without violating the tensor product structure. These systems can hence be solved by established PDE solvers. We remark that our derivation of these equations remains more or less formal, and the existence of weak solutions (in the continuous setting) will not be studied. Let us also mention the work~\cite{Kusch2023} which is related to our work insofar as it shows that for hyperbolic problems a continuous formulation of the projector splitting integrator prior to discretization has favorable numerical stability properties.

We then proceed by deriving corresponding discrete equations that allow to solve the system numerically. Specifically, in our implementation we will use stabilized finite elements~\cite{ern_2021} for discretization. We then apply the method to solve the weak Landau damping on a periodic domain (in order to verify the discretization) as well as a linear equation involving a piecewise linear boundary. Here, in addition to the standard projector-splitting approach we will also consider a weak formulation of the unconventional low-rank integrator proposed in~\cite{CerutiLubich2022}. It consists in modifications regarding the update strategy for the low-rank factors and allows to make the whole scheme rank-adaptive~\cite{CerutiKuschLubich2022}, which is important in practice, based on subspace augmentation. Similar strategies are considered in~\cite{Guo_2022} and~\cite{hauck2022}.

Several recent works have been focussing on the conservation properties in dynamical low-rank integrators~\cite{Einkemmer2019,EINKEMMER2020109063,EINKEMMER2021110495,Einkemmer2022a,GuoWei2022,guo2022}, such as mass, momentum and energy, based on modified Galerkin conditions. This important aspect is not yet addressed in the present paper and will remain for future work.

The paper is organized as follows. In section~\ref{sec:dlra} we detail our idea of a weak formulation of the projector-splitting integrator for the Vlasov-Poisson equation which is capable of handling the inflow boundary conditions. From this, in section~\ref{sec:discrete} we obtain discrete equations by restricting to finite-dimensional spaces according to the Galerkin principle, resulting in Algorithm~\ref{alg:lie-trotter}. The rank-adaptive unconventional integrator is also considered (Algorithm~\ref{alg:rauc}). In section~\ref{sec:experiments} results of numerical experiments are presented to illustrate the principle feasibility of our algorithms.

\section{Weak formulation of the projector splitting integrator \label{sec:dlra}}

Our goal is to develop a projector splitting scheme for DLRA of the transport problem~\eqref{eq:transport} that is capable of handling the inflow boundary conditions~\eqref{eq:boundary_condition}. To achieve this, we start out from a weak formulation of~\eqref{eq:transport} including a weak formulation of the boundary condition according to \cite[Sec.~76.3]{ern_2021}. Assume for now that the electrical field is fixed. Let $W$ be an appropriate closed subspace of $L_2(\Omega)$ (e.g.~$W \subseteq H^1(\Omega)$ is sufficient) and denote by $P_W$ the $L_2$-orthogonal projector onto $W$. Then the goal in the weak formulation is to find $f \in C^1([0,T],W)$ with $f(0) = P_{W} f_0$ and such that at every $t \in (0,T)$ it holds
\begin{equation} \label{eq:weak}
\int_\Omega \partial_t f(t) w \, \mathrm d \bm x \mathrm d \bm v  + a(t, f(t),w) = \ell(w)
\textnormal{ for all } w \in W,
\end{equation}
where
\begin{align*}
a(t,u,w) &= \int_\Omega \bm v \cdot \nabla_{\bm x} u  \, w - \bm E (t,\bm x) \cdot \nabla_{\bm v} u \, w \, \mathrm d\bm x \mathrm d \bm v - \int_{\Gamma^-} \bm v \cdot \bm n_x \, u \, w \, \mathrm d s, \\
\ell(w) &= -\int_{\Gamma^-} \bm v \cdot \bm n_x  \, g \, w \, \mathrm d s.
\end{align*}
Here we have taken into account that $\Omega_v = \mathbb R^d$ is  unbounded and hence the normal vectors for $\Omega$ read $\bm n = (\bm n_x, \bm 0)$. While our derivations are based on the above formulation with fixed electric field $\bm E$, the case where $\bm E$ depends on $f$ will be discussed in section~\ref{sub:efield}. The variational formulation~\eqref{eq:weak} has the advantage that the boundary conditions are integrated in the bilinear form and right hand side. It therefore can be easily combined with the so-called time-dependent variational principle, also called Dirac--Frenkel principle, underlying DLRA.

In order to benefit from the separation of variables in DLRA, it will be necessary to impose some restrictions on the boundary and the function $g$ in the inflow condition~\eqref{eq:boundary_condition}. Specifically, we assume that the spatial domain has a piecewise linear boundary, so that one has a decomposition 
\begin{equation} \label{eq:boundary_decomposition}
\Gamma_x := \partial \Omega_x = \bigcup_\nu \Gamma^{(\nu)}_x\, ,
\end{equation}
where each part $\Gamma^{(\nu)}_x$ has a \emph{constant} outer normal vector $\bm n^{(\nu)}_x$. In addition, we require that $g$ itself admits a separation into space and velocity variables by a finite sum of tensor products,
\begin{equation} \label{eq:tensor_g}
    g(t,\bm x, \bm v) = \sum_\mu g_x^{(\mu)}(t, \bm x) \cdot g_v^{(\mu)}(t, \bm v),
\end{equation}
or can at least be approximated in this form.

In the following sections we derive the projector splitting approach of DLRA in the weak formulation and derive effective equations for the low-rank factors in~\eqref{eq:dlra} in the form of Friedrichs' systems (systems of hyperbolic equations in weak formulation) that are amenable to established PDE solvers. In this way the boundary conditions will be incorporated without sacrificing the tensor product structure.

\subsection{DLRA and projector splitting}

We first present the basic idea of the projector splitting approach for the dynamical low-rank solution of the weak formulation~\eqref{eq:weak}. As required for DLRA, we consider a possibly infinite-dimensional tensor product subspace
\[
W = W_x \otimes W_v \subseteq L_2(\Omega_x) \otimes L_2(\Omega_v) = L_2(\Omega)
\]
where we assume that $W_x \subseteq L_2(\Omega_x)$ and $W_v \subseteq L_2(\Omega_v)$ are suitable subspaces to ensure that~\eqref{eq:weak} is well-defined on $W$. For example, $W_x$ and $W_v$ could be subspaces of $H^1$, which would correspond to a mixed regularity with respect to space and velocity.

The corresponding manifold $\mathcal M_r$ of low-rank functions in $W$ is
\begin{equation}\label{eq:Mr}
\mathcal M_{r} = \Big\{ \varphi(\bm x, \bm v) = \sum_{i=1}^r \sum_{j=1}^r X_i(\bm x) S_{ij}V_j(\bm v) \, \Big\vert \, X_i \in W_x, \, V_j \in W_v, \, S_{ij} \in \mathbb R \Big\},
\end{equation}
with the additional conditions that the $X_1,\dots,X_r$ and $V_1,\dots,V_r$ are orthonormal systems in $W_x$ and $W_v$, respectively, and the matrix $S = [S_{ij}]$ has rank $r$. Following the Dirac-Frenkel variational principle, a dynamical low-rank approximation for the weak formulation~\eqref{eq:weak} asks for a function $f_r \in C^1([0,T],W)$ such that $f_r(t) \in \mathcal M_r$ and
\begin{equation} \label{eq:weak_dlra}
\int_\Omega \partial_t f_r(t) w \, \mathrm d \bm x \mathrm d \bm v  + a(t, f_r(t),w) = \ell(w)
\textnormal{\quad for all } w \in T_{f_r(t)} \mathcal M_r,
\end{equation}
for all $t \in (0,T)$, where $T_{f_r(t)} \mathcal M_r$ is the tangent space of $\mathcal M_r$ at $f_r(t)$ specified below. For the initial value one may take $f_r(0) = \tilde f_0$, where $\tilde f_0$ is a rank-$r$ approximation of $f_0$.

The DLRA formulation~\eqref{eq:weak_dlra} can be interpreted as a projection of the time derivative of the solution onto the tangent space, which implicitly restricts it to the manifold. The projector splitting integrator from~\cite{Lubich2014} is now based on the rather peculiar fact that the tangent spaces of $\mathcal M_r$ can be decomposed into two smaller subspaces corresponding to variations in the $X$ and $V$ component only. Specifically, let
\[
f_r(t,\bm x, \bm v) = \sum_{i=1}^r \sum_{j=1}^r X_i(t,\bm x) S_{ij}(t) V_j(t,\bm v)
\]
be in $\mathcal M_r$ at time $t$, then it holds that
\[
T_{f_r(t)} \mathcal M_r =  T_{V(t)} + T_{X(t)},
\]
where 
\[ T_{V(t)} = \Big\{ \sum_{i=1}^r K_i(\bm x) V_i(t,\bm v)\, \Big\vert \, K_i \in W_x \Big\}, \quad T_{X(t)} = \Big\{ \sum_{i=1}^r X_i(t,\bm x) L_i(\bm v)\, \Big\vert \, L_i \in W_v\Big\}.\]
Note that these two subspaces are not complementary as they intersect in the space
\[
T_{X(t),V(t)} = T_{X(t)} \cap T_{V(t)} = \Big\{ \sum_{i=1}^r \sum_{j=1}^r X_i(t,\bm x) \tilde S_{ij} V_j(t,\bm v) \,\Big|\, \tilde S_{ij} \in \mathbb R \Big\}.
\]
Correspondingly, the $L_2$-orthogonal projection onto the full tangent space $T_{f_r(t)} \mathcal M_r$ can be decomposed as 
\begin{equation} \label{eq:projectors}
    P_{f_r(t)} = P_{V(t)} - P_{X(t),V(t)} + P_{X(t)},
\end{equation}
where $P_{V(t)}$, $P_{X(t),V(t)}$, and $P_{X(t)}$ denote the $L_2$-orthogonal projections onto the subspaces $T_{X(t)}$, $T_{X(t),V(t)}$, and $T_{V(t)}$, respectively. 

The projector-splitting integrator performs the time integration of the system according to the decomposition~\eqref{eq:projectors} of the tangent space projector. Sticking to the weak formulation~\eqref{eq:weak_dlra}, one time step from some point $t_0$ with
\begin{equation} \label{eq:f0}
f_r(t_0) = \sum_{i=1}^r \sum_{j=1}^r X_i^{0}(\bm x) S_{ij}^{0} V_j^{0}(\bm v) \in \mathcal M_r
\end{equation}
to a point $t_1=t_0 + \Delta t$ then consists of the following three steps:  
\begin{enumerate}
\item \label{enum:stepK} Solve the system \eqref{eq:weak_dlra} restricted to the subspace $T_{V^0}$ on the time interval $[t_0,t_1]$ with initial condition $f_r(t_0)$. At time $t_1$ one obtains
\[ \hat f(\bm x, \bm v) = \sum_{j=1}^r K_j(\bm x) V_j^0(\bm v) \in T_{V^0}.\]
Find an orthonormal system $X_1^1, \ldots, X_r^1$ for $K_1, \ldots, K_r$ such that 
\[ K_j = \sum_{i=1}^r X_i^1 \hat S_{ij}.\]

\item \label{enum:stepS} Solve the system \eqref{eq:weak_dlra} restricted to the subspace $T_{X^1,V^0}$ on the time interval $[t_0,t_1]$ with initial condition $\hat f$, and taking into account the minus sign in the projector. At time $t_1$ one obtains 
\[ 
\tilde f(\bm x, \bm v) = \sum_{i=1}^r\sum_{j=1}^r X_i^1(\bm x) \tilde S_{ij} V_j^0(\bm v) \in T_{X^1,V^0}.
\]
This step is often interpreted as a backward in time integration step, which is possible if there is no explicit time dependence in the coefficients or boundary conditions. However, in our case the inflow $g$ is in general time dependent which then does not admit for such an interpretation in general.

\item \label{enum:stepL} Solve the system \eqref{eq:weak_dlra} restricted to the subspace $T_{X^1}$ on the time interval $[t_0,t_1]$ with initial condition $\tilde f$. At time $t_1$ one obtains
\[ 
\bar f(\bm x, \bm v) = \sum_{i=1}^r X_i^1(\bm x) L_i(\bm v) \in T_{X^1}.
\]
Find an orthonormal system $V_1^1, \ldots, V_r^1$ for $L_1, \ldots, L_r$ such that 
\[ L_i = \sum_{j=1}^r S_{ij}^1 V_j^1.\]
\end{enumerate}
As the final solution at time point $t_1$ one then takes
\[ 
f_r(t_1, \bm x, \bm v) \approx \bar f(\bm x, \bm v) = \sum_{i=1}^r \sum_{j=1}^r X_i^1(\bm x) S_{ij}^1 V_j^1(\bm v).
\]

A modification of this scheme is the so called unconventional integrator proposed in~\cite{CerutiLubich2022}. There, the $K$-step and $L$-step are performed independently using the same initial data from~$f_r(t_0)$ (i.e.~the $K$-step is identical with the one above, but the $L$-step is performed on $T_{X^0}$ with initial value~$f_r(t_0)$). As a result, one obtains two new orthonormal sets of component functions $X_1^1,\dots,X_r^1$ and $V_1^1,\dots,V_r^1$ for the spatial and velocity domains. The $S$-step is then performed afterwards in a ``forward'' way (i.e.~without flipping signs) using these new bases. Compared to the standard scheme, this decoupling the $S$-step from the $K$- and $L$-steps also offers a somewhat more natural way of making the scheme rank-adaptive, which is important in practice since a good guess of $r$ might not be known. Specifically, in~\cite{CerutiKuschLubich2022} it is proposed to first augment the new bases to $\hat X = \{X^0_1,\dots,X^0_r,K^1_1,\dots,K^1_r \}$ and $\hat V = \{V^0_1,\dots,V^0_r,L^1_1,\dots,L^1_r \}$ (i.e. including the old bases), orthonormalize these augmented bases, and perform the $S$-step in the space
\[
T_{\hat X, \hat V} = \spann \hat X \otimes \spann \hat V.
\]
Note that this space contains functions of rank $2r$ (in general) and still contains $f_r(t_0)$ which one takes as initial condition for the $S$-step. The solution at time point $t_1$ can then be truncated back to rank $r$ or any other rank (less or equal to $2r$) depending on a chosen truncation threshold using singular value decomposition. 

Regardless of whether the projector splitting or the modified unconventional integrator is chosen, the efficient realization of the above steps is based on their formulation in terms of the component functions $K$, $S$ and $L$ only, which in turn requires a certain separability of space and velocity variables in the system~\eqref{eq:weak_dlra}. We therefore next investigate the single steps in detail, focusing only on the projector splitting integrator and omitting the required modifications for the (rank-adaptive) unconventional integrator (it will be discussed in the discrete setting in section~\ref{sec: rauc}). A particular focus is how to ensure the required separability for the boundary value terms, which will lead to the assumptions~\eqref{eq:tensor_g} and~\eqref{eq:boundary_decomposition} that $g$ is separable and the boundary is piecewise linear. As we will see, this gives effective equations in the form of Friedrichs' systems for the factors $K$ and $L$, and a matrix ODE for $S$. Their discrete counter-parts as well as the algorithmic schemes for both the projector splitting integrator and the unconventional integrator will be presented in section~\ref{sec:discrete}.

\subsection{Weak formulation of the \texorpdfstring{$K$}{K}-step \label{sec:weakK}}

The time dependent solution of the first step, i.e.~\eqref{eq:weak_dlra} restricted to $T_{V^0}$, is a function 
\[
\hat f(t,\bm x, \bm v) = \sum_{j=1}^r K_j(t,\bm x) \cdot V_j^0(\bm v)
\]
on the time interval $[t_0,t_1]$ with initial condition $\hat f(t_0,\bm x,\bm v) = f_r(t_0,\bm x,\bm v)$ from~\eqref{eq:f0}. Here the $V_j^0$ are fixed and the functions $K_1,\dots,K_r$ need to be determined. Their dynamics are governed by a system of first order partial differential equations which we derive in the following.

The test functions $w \in T_{V^0}$ have the form
\[ 
w(\bm x, \bm v) = \sum_{j=1}^r \psi_j^{(x)}(\bm x)\cdot V_j^0(\bm v), 
\qquad \psi_j^{(x)} \in W_x,
\]
for all $j=1,\ldots, r$. Testing with $w$ in~\eqref{eq:weak_dlra} specifically means
\begin{equation} \label{eq:weak_K}
\begin{aligned}
\big(\partial_t \hat f(t) + \bm v \cdot \nabla_{\bm x} \hat f(t) - &\bm E(t,\bm x) \cdot \nabla_{\bm v} \hat f(t), w \big)_{L_2(\Omega, \mathbb R)} \\
&- \big(\bm n_x \cdot \bm v \, \hat f(t), w \big)_{L_2(\Gamma^-, \mathbb R)} = 
-\big(\bm n_x \cdot \bm v \, g, w \big)_{L_2(\Gamma^-, \mathbb R)}.
\end{aligned}
\end{equation}
We will write this equation in vector form and therefore define
\[
\bm K = \big[ K_j \big]_{j=1,\ldots,r}, \qquad \bm \psi^{(x)} = \big[ \psi_j^{(x)} \big]_{j=1,\ldots,r},
\]
which can be regarded as vector-valued functions in $(W_x)^r$.

We investigate the terms in \eqref{eq:weak_K} separately. The term regarding the inner product in $L_2(\Omega, \mathbb R)$ has three parts. The first involves the time derivative:
\begin{align*}
\big(\partial_t \hat f(t) , w\big)_{L_2(\Omega, \mathbb R)}
& = \sum_{i=1}^r \sum_{j=1}^r \int_{\Omega_x} \!\! \partial_t K_i(t,\bm x) \cdot \psi_j^{(x)}(\bm x) \, \mathrm d \bm x \cdot
	\int_{\Omega_v} \! V_i^0(\bm v) \cdot V_j^0(\bm v) \, \mathrm d \bm v \\
& =  \sum_{j=1}^r \int_{\Omega_x} \!\! \partial_t K_j(t,\bm x) \cdot \psi_j^{(x)}(\bm x) \, \mathrm d \bm x \\
&= \big( \partial_t \bm K(t), \bm \psi^{(x)}\big)_{L_2(\Omega_x, \mathbb R^r)},
\end{align*}
where we have used the pairwise orthonormality of $V_1, \dots, V_r$. The second part reads
\begin{align*}
\int_{\Omega} \bm v \cdot \nabla_{\bm x} & \hat f(t,\bm x, \bm v) \cdot w(\bm x, \bm v) \, \mathrm d(\bm x, \bm v) \\ &= \sum_{k=1}^d \sum_{i=1}^r \sum_{j=1}^r \int_{\Omega_x} \!\! \partial_{x_k} K_i(t,\bm x) \cdot \psi_j^{(x)}(\bm x) \, \mathrm d \bm x \cdot \!
	\int_{\Omega_v} \!\! v_k \cdot V_i^0(\bm v) \cdot V_j^0(\bm v) \, \mathrm d \bm v \\
&= \sum_{k=1}^d \big( \mathcal A^k_x \cdot \partial_{x_k} \bm K(t), \bm \psi^{(x)}\big)_{L_2(\Omega_x, \mathbb R^r)}
\end{align*}
with the symmetric $r \times r$ matrix
\[ \mathcal A^k_x = \Big[ \int_{\Omega_v} \! v_k \cdot V_j^0(\bm v) \cdot V_i^0(\bm v) \, \mathrm d \bm v \Big]_{i,j=1,\ldots,r}. \]
Finally, the part involving the electrical field can be written as  
\begin{align*}
\int_{\Omega} & \bm E(t,\bm x) \cdot \nabla_{\bm v} \hat f(t,\bm x, \bm v) \cdot w(\bm x, \bm v) \, \mathrm d(\bm x, \bm v) \\
&= \sum_{i=1}^r \sum_{j=1}^r \int_{\Omega_x}  
	\sum_{k=1}^d \Big[ E_k(t,\bm x) \cdot \int_{\Omega_v} \! \partial_{v_k} V_i^0(\bm v) \cdot V_j^0(\bm v) \, \mathrm d \bm v\Big] 
		\,  K_i(t,\bm x) \cdot \psi_j^{(x)}(\bm x) \, \mathrm d \bm x \\
&= \big( \mathcal K_x(t, \cdot) \cdot \bm K(t), \bm \psi^{(x)}\big)_{L_2(\Omega_x, \mathbb R^r)}
\end{align*}
where
\begin{equation} \label{eq:kappax}
\mathcal K_x(t,\bm x) = \Big[\sum_{k=1}^d E_k(t,\bm x)\,  \int_{\Omega_v} \! \partial_{v_k} V_j^0(\bm v) \cdot V_i^0(\bm v) \, \mathrm d \bm v \Big]_{i,j=1,\ldots,r}\, .
\end{equation}

For the terms in \eqref{eq:weak_K} involving the inflow boundary $\Gamma^-$ we make use of our assumption that the spatial boundary $\partial \Omega_x$ can be decomposed as in \eqref{eq:boundary_decomposition} with a constant normal vector $\bm n_x^{(\nu)}$ on each part $\Gamma_x^{(\nu)}$ of the boundary. We decompose $\Gamma^-$ accordingly into
\[\Gamma^- = \bigcup_\nu \Gamma^{(\nu)}_x \times \Omega^{(\nu)}_v, \qquad
\Omega_v^{(\nu)} = \{\bm v \in \Omega_v \,|\,  \bm n^{(\nu)}_x \cdot \bm v < 0\} .\]
Then the boundary term from the bilinear form can be written as
\begin{align*}
    \int_{\Gamma^-} \! \bm n_x \cdot & \bm v \, \hat f(t, \bm x, \bm v)) \cdot w(\bm x, \bm v) \, \mathrm ds \\
    &= \sum_{i=1}^r\sum_{j=1}^r \sum_{\nu} \int_{\Gamma^{(\nu)}_x} 
    	\Big[ \int_{\Omega^{(\nu)}_v} \! \bm n^{(\nu)}_x \cdot \bm v \, V_i^0(\bm v) \, V_j^0(\bm v) \, \mathrm d \bm v \Big] 
    	\, K_i(t, \bm x) \cdot \psi_j^{(x)}(\bm x) \, \mathrm d s \\
    &= \big(\mathcal B_x(\cdot) \, \bm K(t), \bm \psi^{(x)}\big)_{L_2(\Gamma_x, \mathbb R^r)}
\end{align*}
where
\[ 	\mathcal B_x(\bm x) = 
		\Big[ \sum_\nu \chi_{\Gamma_x^{(\nu)}}(\bm x) \int_{\Omega^{(\nu)}_v} \bm n^{(\nu)}_x \cdot \bm v \, V_j^0(\bm v) \, V_i^0(\bm v) \, \mathrm d \bm v \Big]_{i,j=1,\ldots,r}
\]
with characteristic function $\chi$.

For the boundary term on the right hand side we recall the decomposition \eqref{eq:tensor_g} of the function~$g$ into a sum of tensor products. Proceeding as above we obtain
\begin{align*}
    \int_{\Gamma^-} \!\! \bm n_x \cdot & \bm v \, g(t, \bm x, \bm v) \cdot w(\bm x, \bm v) \, \mathrm ds \\
    &= \sum_{j=1}^r \sum_{\nu, \mu} \int_{\Gamma^{(\nu)}_x} 
    	\Big[ \int_{\Omega^{(\nu)}_v} \! \bm n^{(\nu)}_x \cdot \bm v \, g^{(\mu)}_v(t, \bm v) \, V_j^0(\bm v) \, \mathrm d \bm v \Big] 
    	\, g^{(\mu)}_x(t, \bm x) \cdot \psi_j^{(x)}(\bm x) \, \mathrm d s \\
    &= \big(\bm G_x(t,\cdot), \bm \psi^{(x)}\big)_{L_2(\Gamma_x, \mathbb R^r)}
\end{align*}
with
\[ \bm G_x(t, \bm x) = 
		\Big[ \sum_{\mu,\nu} \chi_{\Gamma_x^{(\nu)}}(\bm x)\, g^{(\mu)}_x(t, \bm x)
			\int_{\Omega^{(\nu)}_v} \! \bm n^{(\nu)}_x \cdot \bm v \, \, g^{(\mu)}_v(\bm v) \,  V_j^0(\bm v) \, \mathrm d \bm v \Big]_{j=1,\ldots,r}\, .
\]

In summary \eqref{eq:weak_K} is equivalent to the following system of first order partial differential equations for $\bm K(t)$ in weak formulation with boundary penalty,
\begin{equation} \label{eq:weakK}
\begin{aligned}
&\big( \partial_t  \bm K(t) + \sum_{k=1}^d \mathcal A_x^k \cdot \partial_k \bm K(t) - \mathcal K_x(t,\cdot) \, \bm K(t), \bm \psi^{(x)}\big)_{L_2(\Omega_x, \mathbb R^r)} \\
&- \big(\mathcal B_x(\cdot) \, \bm K, \bm \psi^{(x)}\big)_{L_2(\Gamma_x, \mathbb R^r)} 
 = -\big(\bm G_x(t,\cdot), \bm \psi^{(x)}\big)_{L_2(\Gamma_x, \mathbb R^r)} \quad \textnormal{for all } \bm \psi^{(x)} \in (W_x)^r,
\end{aligned}
\end{equation}
which needs to be solved for $\bm K(t_1)$ with the initial values
\[
\bm K(t_0) = \Big[ \sum_{i=1}^r X_i^0(\bm x) S_{ij}^0 \Big]_{ j=1,\dots,r}
\]
from the previous approximation~\eqref{eq:f0}.

\subsection{Formulation of the \texorpdfstring{$S$}{S}-step \label{sub:weakS}}

The solution of the second step of the splitting integrator is a time-dependent function 
\[
\tilde f(t, \bm x, \bm v) = \sum_{i=1}^r\sum_{j=1}^r X_i^1(\bm x) S_{ij}(t) V_j^0(\bm v)
\]
in the subspace space $T_{X^1,V^0}$, where the evolution for $S$ on the interval $[t_0, t_1]$ is governed by~\eqref{eq:weak_dlra} restricted to test functions $w \in T_{X^1, V^0}$, but taking into account the minus sign of $P_{X^1,V^0}$ in the corresponding projector splitting \eqref{eq:projectors}. The initial condition reads $\tilde f(t_0,\bm x, \bm v) = \hat f(t_1,\bm x, \bm v) = \sum_{i,j = 1}^r X_i^1(\bm x) \hat S_{ij} V_j^0(\bm v)$. The test function in $w \in T_{X^1,V^0}$ will be written as
\[
w(\bm x, \bm v) = \sum_{m=1}^r \sum_{n=1}^r X_m^1(\bm x) \Sigma_{mn}^{} V_n^0(\bm v), \qquad \Sigma = [\Sigma_{mn}] \in \mathbb R^{r \times r}.
\]
For such test functions, we again consider the different contributions as in~\eqref{eq:weak_K} (with $\tilde f$ instead of~$\hat f$). For the time derivative we get
\begin{equation*} 
\big(\partial_t \tilde f(t) , w\big)_{L_2(\Omega, \mathbb R)} = \sum_{ijmn} \dot S_{ij} \Sigma_{mn} \int_{\Omega_x} X_i^1(\bm x) X_m^1(\bm x) \, \mathrm d \bm x 
	\cdot \int_{\Omega_v} V_j^0(\bm v) V_n^0 (\bm v)\, \mathrm d \bm v = \big(\dot S, \Sigma \big)_F,
\end{equation*}
where $(\cdot,\cdot)_F$ is the Frobenius inner product on $\mathbb R^{r\times r}$.
The first two contributions of the bilinear form read 
\begin{align*}
\int_{\Omega} \bm v \cdot &\nabla_{\bm x} \tilde f(t,\bm x, \bm v) \cdot w(\bm x, \bm v) \, \mathrm d(\bm x, \bm v) \\ &= \sum_{k=1}^d \sum_{ijmn} 
	\underbrace{\int_{\Omega_x} \partial_{x_k} X_i^1(\bm x) \cdot X_m^1(\bm x) \, \mathrm d \bm x}_{{}\coloneqq [D^{(k,1)}]_{mi}}  \cdot
	\underbrace{\int_{\mathbb R^d} v_k \cdot V_j^0(\bm v) \cdot V_n^0(\bm v) \, \mathrm d \bm v}_{{}\coloneqq[C^{(k,1)}]_{nj}} 
	S_{ij} \Sigma_{mn}\\
& = \Big( \sum_{k=1}^d  D^{(k,1)}S (C^{(k,1)})^T, \Sigma\Big)_F
\end{align*}
and
\begin{align*}
\int_{\Omega} & \bm E(t,\bm x) \cdot \nabla_{\bm v} f(t,\bm x, \bm v) \cdot w(\bm x, \bm v) \, \mathrm d(\bm x, \bm v) \\
&= \sum_{k=1}^d \sum_{ijmn} 
	\underbrace{\int_{\Omega_x} E_k(t,\bm x) X_i^1(\bm x) \cdot X_m^1(\bm x) \, \mathrm d \bm x}_{{}\coloneqq [D^{(k,2)}(t)]_{mi}}  \cdot
	\underbrace{\int_{\mathbb R^d} \partial_{v_k} \cdot V_j^0(\bm v) \cdot V_n^0(\bm v) \, \mathrm d \bm v}_{{}\coloneqq [C^{(k,2)}]_{nj}} 
	S_{ij} \Sigma_{mn}\\
&= \Big( \sum_{k=1}^d  D^{(k,2)}S (C^{(k,2)})^T, \Sigma\Big)_F.
\end{align*}

The terms stemming from the boundary condition take the form
\begin{align*}
    \int_{\Gamma^-} \! \bm n_x &\cdot \bm v \, f(t, \bm x, \bm v) \cdot w(\bm x, \bm v) \, \mathrm ds \\
    &= \sum_{\nu} \sum_{ijmn}  
    	\underbrace{\int_{\Gamma^{(\nu)}_x} X_i^1(\bm x) \cdot X_m^1(\bm x) \, \mathrm d s}_{{}\coloneqq [B^{(\nu,x)}]_{mi}} \cdot 
    	\underbrace{\int_{\Omega^{(\nu)}_v} \bm n^{(\nu)}_x \cdot \bm v \, V_j^0(\bm v) \, V_n^0(\bm v) \, \mathrm d \bm v}_{{}\coloneqq [B^{(\nu,v)}]_{nj}} \,
    	S_{ij} \Sigma_{mn} \\
    &= \sum_{\nu} \big(  B^{(\nu,x)}S (B^{(\nu,v)})^T, \Sigma\big)_F \, ,
\end{align*}
and
\begin{align*}
    \int_{\Gamma^-} \! &\bm n_x \cdot \bm v \, g(t, \bm x, \bm v) \cdot w(\bm x, \bm v) \, \mathrm ds \\
    &= \sum_{\nu \mu} \sum_{mn} 
    \underbrace{\int_{\Gamma^{(\nu)}_x} g^{(\mu)}_x(t, \bm x) \cdot X_m^1(\bm x) \, \mathrm d s}_{{}\coloneqq [\bm g^{(\nu,\mu,x)}(t)]_m} \cdot
    \underbrace{\int_{\Omega^{(\nu)}_v} \bm n^{(\nu)}_x \cdot \bm v \, g^{(\mu)}_v(\bm v) \, V_n^0(\bm v) \, \mathrm d \bm v}_{{}\coloneqq [\bm g^{(\nu,\mu,v)}(t)]_n}\,
    	\Sigma_{mn} \\
    &= \big(G_S ,\Sigma \big)_F, \quad
\end{align*}
with 
\[
G_S = \sum_{\nu \mu} \bm g^{(\nu,\mu,x)}\cdot\big(\bm g^{(\nu,\mu,v)}\big)^T \,.
\]

Putting everything together and testing with all $\Sigma = [\Sigma_{mn}] \in \mathbb R^{r \times r}$ yields the ODE
\begin{equation} \label{eq:weakS}
\dot S - \sum_{k=1}^d\Big[  D^{(k,1)}S (C^{(k,1)})^T - D^{(k,2)}S (C^{(k,2)})^T\Big]
    + \sum_\nu  B^{(\nu,x)}S (B^{(\nu,v)})^T  =  G_S
\end{equation} 
(in strong form) for obtaining $S(t_1) = \tilde S$ from the initial condition $S(t_0) = \hat S$. Note again that compared to the original problem~\eqref{eq:weak_dlra} all signs (except the one of $\dot S$) have been flipped due to the negative sign of the corresponding projector in the splitting~\eqref{eq:projectors}.

\subsection{Weak formulation of the \texorpdfstring{$L$}{L}-step}\label{sec:weakL}

The equations for the $L$-step are derived in a similar way as the $K$-step in section~\ref{sec:weakK}. Here we seek
\[ 
\bar f(t, \bm x, \bm v) = \sum_{i=1}^r X_i^1(\bm x) \cdot L_i(t,\bm v) 
\]
in the subspace $T_{X^1}$ with fixed $X_1^1,\dots,X_r^1$ from the initial value $\bar f(t_0,\bm x,\bm v) = \tilde f(t_1,\bm x, \bm v)$. The test functions in $T_{X^1}$ are of the form
\[
w(\bm x, \bm v) = \sum_{i=1}^r X_i^1(\bm x) \cdot \psi_i^{(v)}(\bm v), \qquad
\psi_i^{(v)} \in W_v. 
\]
Let again
\[
\bm L = \big[ L_i \big]_{i=1,\ldots,r}, \qquad \bm \psi^{(v)} = \big[ \psi_i^{(v)} \big]_{i=1,\ldots,r}\, .
\]
be vector valued functions in $(W_v)^r$.

When restricting the weak formulation as in~\eqref{eq:weak_K} (with $\bar f$ instead of~$\hat f$) to test functions in~$T_{X^1}$, the first three resulting contributions take a similar form as in the $K$-step,
\[
\big(\partial_t \bm L(t) - \sum_{k=1}^d \mathcal A^k_v(t) \cdot \partial_k \bm L(t)
+ \mathcal K_v(\cdot) \cdot \bm L(t), \bm \psi^{(v)}\big)_{L_2(\Omega_v, \mathbb R^r)}
\]
but this time with
\begin{equation} \label{eq:kappav}
\begin{aligned}
 \mathcal A^k_v(t) &= \Big[ \int_{\Omega_x} \! E_k(t,\bm x) \,  X_j^1(\bm x) \cdot X_i^1(\bm x) \, \mathrm d \bm x \Big]_{i,j=1,\ldots,r}, \\
 \mathcal K_v(\bm v) &= \Big[\sum_{k=1}^d v_k \, \int_{\Omega_x} \! \partial_{x_k} X_j^1(\bm x) \cdot X_i^1(\bm x) \, \mathrm d \bm x \Big]_{i,j=1,\ldots,r}
\end{aligned}
\end{equation}
(note that in the $K$-step $\mathcal A_x^k$ featured the velocity and $\mathcal K_x$ the electric field).

The boundary term in the bilinear form gives
\begin{align*}
    \int_{\Gamma^-} \! \bm n_x \cdot \bm v &\, f(t, \bm x, \bm v)) \cdot w(\bm x, \bm v) \, \mathrm ds \\
    &= \sum_{i=1}^r \sum_{j=1}^r \sum_{\nu} 
    \int_{\Omega^{(\nu)}_v}  \bm n^{(\nu)}_x \cdot \bm v \, 
    	\Big[ \int_{\Gamma^{(\nu)}_x}
    	  \, X_i^1(\bm x) \, X_j^1(\bm x) \, \mathrm d s \Big] 
    	\, L_i(t, \bm v) \cdot \psi_j^{(v)}(\bm v) \, \mathrm d \bm v \\ 
&= \big(\mathcal B_v(\cdot) \, \bm L, \bm \psi^{(v)}\big)_{L_2(\Omega_v, \mathbb R^r)} \, ,
\end{align*}
where  
\[ 	\mathcal B_v(\bm v) = 
		\Big[ \sum_\nu \chi_{\Omega^{(\nu)}_v}(\bm v) \bm n^{(\nu)}_x \cdot \bm v \int_{\Gamma^{(\nu)}_x}
    	  X_j^1(\bm x) \, X_i^1(\bm x) \, \mathrm d s\Big]_{i,j=1,\ldots,r}.
\]
Hence, the boundary condition on the spatial domain $\Omega_x$ leads to the additional multiplicative term on the whole domain $\Omega_v$. Boundary terms for the velocity variable are not present, since since we assume $\Omega_v=\mathbb R^d$ to be unbounded.

For the right hand side, we have
\begin{align*}
    \int_{\Gamma^-} \!\! & \bm n_x \cdot \bm v  \, g(t, \bm x, \bm v) \cdot w(\bm x, \bm v) \, \mathrm ds \\
    &= \sum_{i=1}^r \sum_{\mu\nu} \int_{\Omega_v} \chi_{\Omega^{(\nu)}_v}(\bm v) \bm n^{(\nu)}_x \cdot \bm v \,
    	\Big[ \int_{\Gamma^{(\nu)}_x}  g^{(\mu)}_x(t,\bm x) \, X_i^1(\bm x) \, \mathrm d s \Big] 
    	\, g^{(\mu)}_v(t, \bm v) \cdot \psi_i^{(v)}(\bm v) \, \mathrm d \bm v \\
    &=\big(\bm G_v(t, \cdot), \bm \psi^{(v)}\big)_{L_2(\Omega_v, \mathbb R^r)} 
\end{align*}
with
\[ 
\bm G_v(t, \bm v) = 
		\Big[ \sum_{\mu\nu} g^{(\mu)}_v(t, \bm v) \chi_{\Omega^{(\nu)}_v}(\bm v) \bm n^{(\nu)}_x \cdot \bm v
			\int_{\Gamma^{(\nu)}_x} g^{(\mu)}_x(t,\bm x) \,  X_i^1(\bm x) \, \mathrm d s \Big]_{i=1,\ldots,r}.
\]

The resulting equation for the $L$-step reads: Find $\bm L(t) \in (W_v)^r$ for $t \in [t_0,t_1]$ such that
\begin{equation} \label{eq:weakL}
\begin{aligned}
    \big(\partial_t \bm L(t) &- \sum_{k=1}^d \mathcal A^k_v(t) \cdot \partial_k \bm L(t)
+ \mathcal K_v(\cdot) \, \bm L(t) - \mathcal B_v(\cdot) \bm L, \bm \psi^{(v)}\big)_{L_2(\Omega_v, \mathbb R^r)}\\
&= -\big(\bm G_v(t, \cdot), \bm \psi^{(v)}\big)_{L_2(\Omega_v, \mathbb R^r)}
\quad \textnormal{for all } \bm \psi^{(v)} \in (W_v)^r.
\end{aligned}
\end{equation}
The initial condition is
\[
\bm L(t_0) = \Big[ \sum_{i=1}^r \tilde S_{ij}^0 V^0_j(\bm v) \Big]_{ j=1,\dots,r},
\]
where the $V_j^0$ are from $f_r(t_0)$ in the previous time step~\eqref{eq:f0}, and $\tilde S$ has been determined in the $S$-step.

\subsection{Electrical field \label{sub:efield}}

So far we have assumed that the electrical field does not depend on the density $f$. However, in the setting of a Lie splitting the other case can be included quite easily as well. At the beginning of the time step from $t_0$ to $t_1$, as discussed in sections~\ref{sec:weakK}--\ref{sec:weakL}, the electrical field can be computed via the Poisson equation \eqref{eq:poisson} and then simply held constant during that step. This will only introduce an error of first order in time, as does the projector splitting itself.

\section{Discrete equations \label{sec:discrete}}

In sections~\ref{sec:weakK}--\ref{sec:weakL} the governing equations for the projector splitting integrator have been formulated in a continuous setting. They consist of the two Friedrichs' systems~\eqref{eq:weakK} and~\eqref{eq:weakL} for $\bm K$ and $\bm L$, and the finite-dimensional system of ODEs~\eqref{eq:weakS} for $S$. In this section we will introduce a finite element discretization for the approximate numerical solution of the hyperbolic systems and derive the governing discrete equations that need to be solved in the practical computation.

A natural choice for the discretization is to use the discontinuous Galerkin method. However, the resulting equations for $\bm K$, $\bm L$, and $S$ (Sec.~\ref{sec:dlra}) include derivatives of the bases functions $X_i$ and $V_j$, see for example Eqns.~\eqref{eq:kappav} and \eqref{eq:kappax}. Hence, for a first approach we will use continuous finite element method, which of course has to be stabilized. It remains to future work to investigate the use of DG methods. 

\subsection{Discretization}

Instead of an unbounded domain we will now use a finite domain $\Omega_v = [-v_{\max},v_{\max}]^d$ for the velocity variable with $v_{\max}$ big enough, and impose periodic boundary conditions.

Let meshes $\mathcal T_{x/v}$ on the domains $\Omega_{x/v}$ be given and define $H^1$-conforming discretization spaces
\[ 
W_{x}^h = \mathrm{span}\{ \varphi^{(x)}_\alpha \,|\, \alpha=1,\ldots,n_{x}\}, \qquad W_{v}^h = \mathrm{span}\{ \varphi^{(v)}_\beta \,|\, \beta=1,\ldots,n_v\} \, .
\]
We are now looking for an approximate solution of the Vlasov--Poisson equation~\eqref{eq:transport} in the manifold
\begin{equation*} \label{eq:Mrh}
\mathcal M_r^h = \Big\{ \varphi^h(\bm x, \bm v) = \sum_{i=1}^r \sum_{j=1}^r X_i^h(\bm x) S_{ij} V_j^h(\bm v) \, \Big\vert \, X_i^h \in W_x^h, \, V_j^h \in W_v^h, \, S_{ij} \in \mathbb R \Big\},
\end{equation*}
of rank-$r$ finite element functions, where as before the systems $X_1^h,\ldots,X_r^h$ and $V_1^h,\dots,V_r^h$ are assumed to be orthonormal, and $S = [S_{ij}]$ has rank $r$. To represent elements in $\mathcal M_r^h$, we introduce coefficient matrices
\[
\mathsf X = [\mathsf X_{\alpha i}] \in \mathbb R^{n_x \times r}, \qquad \mathsf V = [\mathsf V_{\beta j}] \in \mathbb R^{n_v \times r}
\]
such that
\begin{equation}\label{eq:discreteXV}
X^h_i(\bm x) = \sum_{\alpha=1}^{n_x} \mathsf X_{\alpha i} \varphi^{(x)}_\alpha(\bm x), \qquad
V^h_j(t,\bm v) = \sum_{\beta=1}^{n_v} \mathsf V_{\beta j} \varphi^{(v)}_\beta(\bm v).
\end{equation}
The goal is then to find the solution in the form
\[
f_r^h(t, \bm x, \bm v) = \sum_{\alpha = 1}^{n_x} \sum_{\beta = 1}^{n_v} \mathsf f_{\alpha \beta}^{}(t) \varphi^{(x)}_\alpha(\bm x) \varphi^{(v)}_\beta(\bm v)
\]
with coefficients
\[
\mathsf f(t) = [\mathsf f_{\alpha \beta}(t)] = \mathsf X(t) \cdot S(t) \cdot \mathsf V(t)^T \in \mathbb R^{n_x \times n_v}
\]
at prescribed time points $t$. Note that in order to have the  $X_1^h,\dots,X_k^h$ orthonormal in $L_2$ as was always assumed above, the matrix $\mathsf X$ needs to be orthogonal with respect to the inner product of the mass matrix of the basis functions $\varphi^{(x)}_\alpha$ (the matrix $\mathsf M_x$ below). A similar remark applies to $\mathsf V$.

Following the Galerkin principle, the steps of the projector splitting integrator can be formulated in terms of the matrices $\mathsf X$, $S$, and $\mathsf V$ by solving the equations~\eqref{eq:weakK},~\eqref{eq:weakS}, and~\eqref{eq:weakL} with the finite element spaces $W_{x/v}^h$ instead of $W_{x/v}$. This means that we use the finite dimensional spaces $(W_{x/v}^h)^r$ both as ansatz and test spaces. Specifically, for the Friedrichs' systems~\eqref{eq:weakK} and~\eqref{eq:weakL} our aim is to find factors $K^h_j$ and $L^h_i$ of the form
\begin{equation}\label{eq:discreteKL}
K^h_j(t,\bm x) = \sum_{\alpha=1}^{n_x} \mathsf K_{\alpha j}(t) \varphi^{(x)}_\alpha(\bm x), \qquad
L^h_i(t,\bm v) = \sum_{\beta=1}^{n_v} \mathsf L_{\beta i}(t) \varphi^{(v)}_i(\bm v)
\end{equation}
for $i,j=1,\dots,r$. We gather the coefficients in the matrices 
\begin{equation*} \label{eq:coefficents}
    \mathsf K(t) = \big[ \mathsf K_{\alpha j}(t)\big]_{\alpha ,j} \in \mathbb R^{n_x \times r}, \quad
\mathsf L(t) = \big[ \mathsf L_{\beta i}(t) \big]_{\beta,i} \in \mathbb R^{n_v \times r} \,.
\end{equation*} 
The governing discrete equations are obtained by inserting the expressions~\eqref{eq:discreteXV} for $X_i$ and $V_j$ as well as~\eqref{eq:discreteKL} for $K_i$ and $L_j$ into~\eqref{eq:weakK},~\eqref{eq:weakS}, and~\eqref{eq:weakL}. In order to formulate them conveniently we will use the following additional discretization matrices:
\begin{align*}
\mathsf M_x &= \big[ \langle \varphi^{(x)}_\alpha, \, \varphi^{(x)}_{\alpha'} \rangle_{L_2(\Omega_x)} \big]_{\alpha, \alpha'}, &
\mathsf M_{x, E_k} &= \big[ \langle \varphi^{(x)}_\alpha, \,  E_k \, \varphi^{(x)}_{\alpha'} \rangle_{L_2(\Omega_x)} \big]_{\alpha, \alpha'}, \\
\mathsf T_{x,k} &= \big[ \langle \varphi^{(x)}_\alpha, \, \partial_{x_k} \varphi^{(x)}_{\alpha'} \rangle_{L_2(\Omega_x)} \big]_{\alpha, \alpha'}, &
\mathsf M_{x,\Gamma^{(\nu)}_x} &= \big[ \langle \varphi^{(x)}_\alpha, \, \varphi^{(x)}_{\alpha'}  \rangle_{L_2(\Gamma^\nu)} \big]_{\alpha, \alpha'}, \\
\mathsf M_{v} & = \big[ \langle \varphi^{(v)}_\beta, \, \varphi^{(v)}_{\beta'} \rangle_{L_2(\Omega_v)} \big]_{\beta,\beta'}, &
\mathsf M_{v,k} &= \big[ \langle \varphi^{(v)}_\beta, \, v_k \varphi^{(v)}_{\beta'} \rangle_{L_2(\Omega_v)} \big]_{\beta,\beta'}, \\
\mathsf T_{v,k} &= \big[ \langle \varphi^{(v)}_\beta, \, \partial_{v_k} \varphi^{(v)}_{\beta'} \rangle_{L_2(\Omega_v)} \big]_{\beta, \beta'}, &
\mathsf M_{v,\Omega^{(\nu)}} &= \big[ \langle \varphi^{(v)}_\beta, \, n^{(\nu)} \cdot v \, \varphi^{(v)}_{\beta'} \rangle_{L_2(\Omega^{(\nu)}_v)} \big]_{\beta,\beta'}, \\
\mathsf G_x(t) &= \big[ \langle \varphi^{(x)}_\alpha , \, G_{x,j}(t,\cdot) \rangle_{L_2(\Gamma_x, \mathbb R^r)} \big]_{\alpha,j}, &
\mathsf G_v(t) &= \big[ \langle \varphi^{(v)}_\beta , \, G_{v,i}(t,\cdot) \rangle_{L_2(\Omega_v, \mathbb R^r)} \big]_{\beta,i}.
\end{align*}

In addition, we will take into account that finite element solutions of hyperbolic partial differential equations have to be stabilized. For that purpose we use the continuous interior penalty (CIP) stabilization \cite{ern_2021}. For a mesh $\mathcal T$, its bilinear form is given by
\begin{equation*}
    s_{\mathcal T}^{\mathrm{CIP}}(v_h, w_h) = 
    \sum_{F \in \mathcal F_0^h} h_F^2 \left( [\![\nabla v_h ]\!], [\![ \nabla w_h ]\!] \right)_{L_2(F)} 
\end{equation*}
where $\mathcal F_0^h$ is the set of all interior interfaces of $\mathcal T$ and $[\![\cdot ]\!]$ is the jump term with respect to the interface $F$. The corresponding discretization matrices are denoted by
\begin{equation*}
    \mathsf C_{x} = \big[ s_{\mathcal T_x}^{\mathrm{CIP}}(\varphi^{(x)}_j, \varphi^{(x)}_i) \big]_{i,j}, \qquad
    \mathsf C_{v} = \big[ s_{\mathcal T_v}^{\mathrm{CIP}}(\varphi^{(v)}_j, \varphi^{(v)}_i) \big]_{i,j} \,.
\end{equation*}

Eventually, the discrete versions of equations~\eqref{eq:weakK},~\eqref{eq:weakS}, and~\eqref{eq:weakL} read as follows:
\begingroup
\allowdisplaybreaks
\begin{subequations}
\begin{align}
    \mathsf M_x \dot{\mathsf K} 
        &= -\sum_{k=1}^d \big( \mathsf T_{x,k} \cdot \mathsf K \cdot \langle \mathsf M_{v,k} \rangle_{\mathsf V}^T
            - \mathsf M_{x,E_k} \cdot \mathsf K \cdot \langle \mathsf T_{v,k} \rangle_{\mathsf V}^T \big) 
            - \delta \, \mathsf C_x \cdot \mathsf K \notag \\
        &\quad + \sum_\nu \mathsf M_{x,\Gamma^{(\nu)}_x} \cdot\mathsf K \cdot \langle \mathsf M_{v,\Omega^{(\nu)}_v}\rangle_{\mathsf V}^T 
                    - \mathsf G_x(t), \label{eq:K}\\
    \dot{S} &= \sum_{k=1}^d \big(  \langle \mathsf T_{x,k} \rangle_{\mathsf X} \cdot S \cdot \langle \mathsf M_{v,k} \rangle_{\mathsf V}^T
        - \langle \mathsf M_{x,E_k} \rangle_{\mathsf X} \cdot S \cdot \langle \mathsf T_{v,k} \rangle_{\mathsf V}^T\big) \notag \\
        &\quad - \sum_\nu \langle \mathsf M_{x,\Gamma^{(\nu)}_x} \rangle_{\mathsf X} \cdot S \cdot \langle \mathsf M_{v,\Omega^{(\nu)}_v}\rangle_{\mathsf V}^T + G_S, \label{eq:S} \\[\medskipamount]
    \mathsf M_v \dot{\mathsf L}
        &= -\sum_{k=1}^d \big( \mathsf M_{v,k} \cdot \mathsf L \cdot  \langle \mathsf T_{x,k} \rangle_{\mathsf X}^T -\mathsf T_{v,k} \cdot \mathsf L \cdot \langle \mathsf M_{x,E_k} \rangle_{\mathsf X}^T \big)
            - \delta \, \mathsf C_v \cdot \mathsf L \notag \\
        &\quad + \sum_\nu \mathsf M_{v,\Omega^{(\nu)}_v} \cdot\mathsf L \cdot \mathsf \langle \mathsf M_{x,\Gamma^{(\nu)}_x}\rangle_{\mathsf X}^T 
                    - \mathsf G_v(t) \label{eq:L} \,.
\end{align}
\end{subequations}
\endgroup
Here we use the abbreviations
\[ \langle \mathsf M\rangle_{\mathsf X} = \mathsf X^T \cdot \mathsf M \cdot \mathsf X, \quad 
\langle \mathsf M\rangle_{\mathsf V} = \mathsf V^T \cdot \mathsf M \cdot \mathsf V\]
for discretization matrices $\mathsf M$. The parameter  $\delta \ge 0$ controls the stabilization.

\subsection{Discrete projector splitting integrator}\label{sec: algorithm}

Based on the above equations, the complete realization of the projector splitting scheme, including the appropriate initial conditions and the orthogonalization steps, is outlined in Algorithm~\ref{alg:lie-trotter}.

It also includes in lines 1 and 2 the computation of the electric field as described in section~\ref{sub:efield}. For that purpose the density $\rho$ has to be calculated given by $f_r^h(t_0,\cdot,\cdot) \in \mathcal M_r^h$ by integrating out the $v$ variable. The result is a finite element function $\rho^h \in W_x^h$ which is used on the right hand side of the Poisson equation~\eqref{eq:poisson}. We used ansatz functions of order two to solve the resulting linear equation directly for the potential $\Phi^h$. The electrical field $\bm E^h$ can be computed as the negative gradient of $\Phi^h$ and is a discontinuous finite element function on the same mesh. Now, updating the discretization matrices $\mathsf M_{x,E_k}$ involving the electrical field, the projector splitting steps from $t_0$ and $t_1$ can be performed.

\begin{algorithm}
\caption{First order dynamical low rank integrator with projector splitting \label{alg:lie-trotter}}
\begin{algorithmic}[1]
\Require $\mathsf X^0, \, S^0, \, \mathsf V^0$
\LineComment update electrical field
\State Solve Poisson equation \eqref{eq:poisson} for density resulting from $\mathsf X^0, S^0, \mathsf V^0$
\State Compute $\bm E$ and update $\mathsf M_{x,E_k}$
\LineComment $K$ step
\State $\mathsf K^0 = \mathsf X^0 S^0$
\State Solve \eqref{eq:K} with $\mathsf V = \mathsf V^0$ on interval $[t_0,t_1]$ with initial condition $\mathsf K^0$ to obtain $\mathsf K^1$
\State Compute orthonormal basis $\mathsf X^1$ with respect to $\mathsf M_x$ such that $\mathsf X^1 \cdot \hat S = \mathsf K^1$
\LineComment $S$ step
\State Solve \eqref{eq:S} with $\mathsf X = \mathsf X^1$, $\mathsf V = \mathsf V^0$ on interval $[t_0,t_1]$ with initial condition $\hat S$ to obtain $\tilde S$
\LineComment $L$ step
\State $\mathsf L^0 = \mathsf V^0 \tilde S^T$
\State Solve \eqref{eq:L} on interval $[t_0,t_1]$ with initial condition $\mathsf L^0$ to obtain $\mathsf L^1$
\State Compute orthonormal basis $\mathsf V^1$  with respect to $\mathsf M_v$ such that $\mathsf V^1 \cdot (S^1)^T = \mathsf L^1$
\State \textbf{return} $\mathsf X^1, \, S^1, \, \mathsf V^1$
\end{algorithmic}
\end{algorithm}

\subsection{Rank adaptive algorithm}\label{sec: rauc}

The rank-adaptive algorithm is based on the unconventional low-rank integrator proposed in~\cite{CerutiLubich2022} and its rank-adaptive extension in~\cite{CerutiKuschLubich2022}. It is displayed in Algorithm~\ref{alg:rauc}. Starting with a rank-$r_0$ function, it first performs independent $K$-steps and $L$-steps from the initial data $\mathsf X^0$ and $\mathsf V^0$. The computed factors $\mathsf K^1$ and $\mathsf L^1$ are used to enlarge the previous bases for $\mathsf X^0$ and $\mathsf V^0$ to dimension $2 r_0$ each. This is achieved by computing orthonormal bases $\hat{\mathsf X}$ and $\hat{\mathsf V}$ of the augmented matrices $[\mathsf X^0, \mathsf K^1] \in \mathbb R^{n_x \times 2 r_0}$ and $[\mathsf V^0, \mathsf L^1] \in \mathbb R^{n_v \times 2 r_0}$. Then a new coefficient matrix of the form
\[
 [\mathsf f_{\alpha \beta}] = \hat{\mathsf X} S \hat{\mathsf V}^T
\]
is sought by performing a `forward' $S$-step
\begin{align}\label{eq:Smodified}
 \dot{S} &= - \sum_{k=1}^d \big(  \langle \mathsf T_{x,k} \rangle_{\hat{\mathsf X}} \cdot S \cdot \langle \mathsf M_{v,k} \rangle_{\hat{\mathsf V}}^T
        - \langle \mathsf M_{x,E_k} \rangle_{\hat{\mathsf X}} \cdot S \cdot \langle \mathsf T_{v,k} \rangle_{\hat{\mathsf V}}^T\big) \notag \\
        &\quad + \sum_\nu \langle \mathsf M_{x,\Gamma^{(\nu)}_x} \rangle_{\hat{\mathsf X}} \cdot S \cdot \langle \mathsf M_{v,\Omega^{(\nu)}_v}\rangle_{\hat{\mathsf V}}^T - G_S, 
\end{align}
which differs from~\eqref{eq:S} in the signs of the right-hand side. For the initial condition $S(t_0) = \hat S$, $\mathsf X^0 S^0 \mathsf V^0$ has to be expressed with respect to the larger bases $\hat{\mathsf X}$ and $\hat{\mathsf V}$, i.e.~as $\hat{\mathsf X} \hat S \hat{\mathsf V}$. The new solution then has rank $2 r_0$ (in general) and can be truncated to a lower rank $r_1$ according to a tolerance, which yield the actual new factors $\mathsf X^1 \in \mathbb R^{n_x \times r_1}$, $S^1 \in \mathbb R^{r_1 \times r_1}$, and $\mathsf V^1 \in \mathbb R^{n_v \times r_1}$.

\begin{algorithm}
\caption{First order rank adaptive unconventional integrator (RAUC) \label{alg:rauc}}
\begin{algorithmic}[1]
\Require $\mathsf X^0, \, S^0, \, \mathsf V^0$, with $S^0 \in \mathbb R^{r_0 \times r_0}$, tolerance $\epsilon$
\LineComment update electrical field
\State Solve Poisson equation \eqref{eq:poisson} for density resulting from $\mathsf X^0, S^0, \mathsf V^0$
\State Compute $\bm E$ and update $\mathsf M_{x,E_k}$
\LineComment Compute augmented $\mathsf X$ basis
\State $\mathsf K^0 = \mathsf X^0 S^0$
\State Solve \eqref{eq:K} on interval $[t_0,t_1]$ with initial condition $\mathsf K^0$ to obtain $\mathsf K^1$
\State Compute orthonormal basis $\hat{\mathsf X}$ with respect to $\mathsf M_x$ such that $\hat{\mathsf X} \cdot [R_x, \tilde R_x] = [\mathsf X^0, \mathsf K^1]$
\LineComment Compute augmented $\mathsf V$ basis
\State $\mathsf L^0 = \mathsf V^0 (S^0)^T$
\State Solve \eqref{eq:L} on interval $[t_0,t_1]$ with initial condition $\mathsf L^0$ to obtain $\mathsf L^1$
\State Compute orthonormal basis $\hat{\mathsf V}$ with respect to $\mathsf M_v$ such that $\hat{\mathsf V} \cdot [R_v, \tilde R_v] = [\mathsf V^0, \mathsf L^1]$
\LineComment $S$ step
\State Compute initial condition $\hat S = R_x \cdot S^0 \cdot R_v^T$
\State Solve modified equation~\eqref{eq:Smodified} on interval $[t_0,t_1]$ with initial condition $\hat S$ to obtain $\tilde S$
\LineComment truncation
\State Compute SVD $Q_x \cdot \Sigma \cdot Q_v^T = \tilde S$\ with $\Sigma = \mathrm{diag}(\sigma_i)$ and monotonically decreasing $\sigma_i$
\State For tolerance $\epsilon$ compute $r_1 \le 2 r_0$ as the minimal number such that
\[ \sum_{i=r_1+1}^{2 r_0} \sigma_i^2 < \epsilon^2. \vspace*{-3ex}\]
\State Set $\mathsf X^1 = \hat{\mathsf X} \cdot Q_x(:,1:r_1)$, $\mathsf V^1 = \hat{\mathsf V} \cdot Q_v(:,1:r_1)$,  $S^1 = \mathrm{diag}(\sigma_i)_{i=1,\ldots,r_1}$
\State \textbf{return} $\mathsf X^1, \, S^1, \, \mathsf V^1$
\end{algorithmic}
\end{algorithm}

\section{Numerical experiments \label{sec:experiments}}

We present numerical results of two experiments for testing the methods described in section~\ref{sec:discrete}. Since to our knowledge there does not exist an analytic solution of the Vlasov--Poisson equation on bounded domains to which our numerical results could be compared, we chose to consider the following two setups. The first experiment will be the classical Landau damping where analytic results on the decay of the electric energy are known. However, the domain is assumed to be periodic so that no boundary conditions are needed. This scenario has been treated in several previous works. Our second experiment will then include boundary conditions on a polygonal spatial domain to actually test our proposed approach for handling these.

In our tests linear finite elements are used and the matrix ODEs (lines~4,~6, and~8 in Alg.~\ref{alg:lie-trotter}, and line~11 in Alg.~\ref{alg:rauc}) are solved using an explicit Runge-Kutta method of third order. The implementation is based on the finite element library \verb|MFEM|~\cite{mfem} and uses its Python wrapper \verb|PyMFEM|.\footnote{\url{https://github.com/mfem/PyMFEM}} The computations are carried out using standard numerical routines from \verb|numpy| and \verb|scipy|. All experiments have been performed on a desktop computer (i9 7900, 128 GB RAM), the code is available online.\footnote{\url{https://github.com/azeiser/dlra-bc}}

\subsection{Landau damping}

As a first numerical test we consider the classical Landau damping in two spatial and velocity dimensions. We use periodic domains $\Omega_x = [0, 4 \pi]^2$ and $\Omega_v = [-6, 6]^2$ and a background density of $\rho_b = 1$. As initial condition we choose
\begin{equation*}
f(0,\bm x,\bm v) = \frac{1}{2\pi} \mathrm e^{-\vert \bm v\vert^2/2}\, \big( 1 + \alpha \cos(k x_1) + \alpha \cos(k x_2)\big), \quad 
\alpha = 10^{-2}, \, k = \frac 1 2.
\end{equation*}
The same setup was investigated in~\cite{einkemmer2018a}. For this case linear analysis shows that the electric field decays with a rate of $\gamma \approx 0.153$. For the discretization we use regular grids with $n_x=64^2$ and $n_v=256^2$ degrees of freedom (level $0$). 

\begin{figure}[t]
\begin{center}
\vspace*{-4ex}
\includegraphics[width=0.95\textwidth]{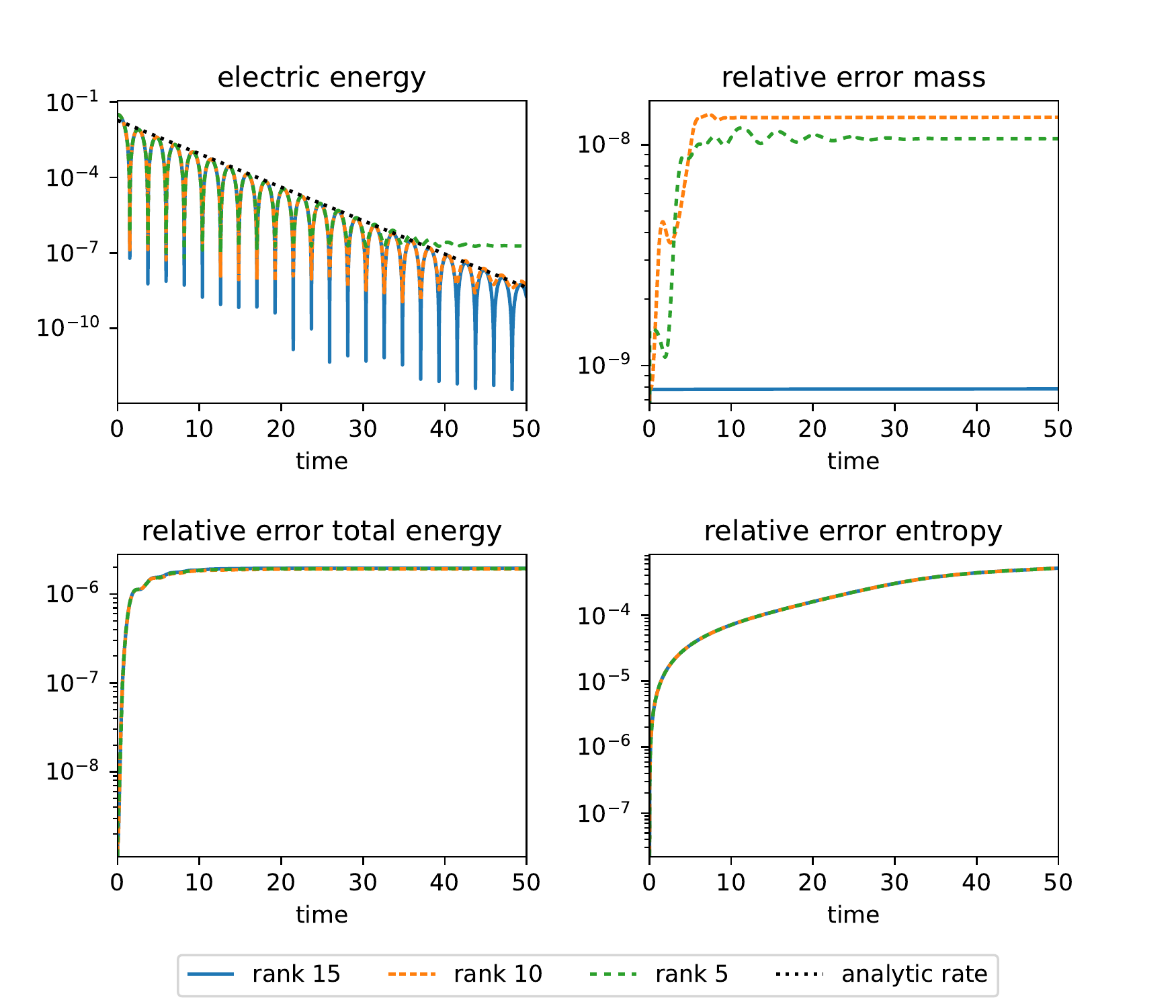}
\end{center}
\caption{Simulation results of the 2+2-dimensional Landau damping using fixed ranks (Alg.~\ref{alg:lie-trotter}); electric energy including the analytical decay rate (upper left) and relative error of the invariants \eqref{eq:invariants}. For the total energy and entropy the lines for all three ranks almost overlap. }
\label{fig:landau_fixed}
\end{figure}

Using Algorithm~\ref{alg:lie-trotter} the simulation was carried out fixed ranks $r=5,10,15$ and a time step of $\Delta t = 0.005$. The results are shown in Figure~\ref{fig:landau_fixed}. As can be seen in the top left plot, for rank $r=5$ the computed electric energy $\frac 1 2 \, \int_{\Omega_x} | \bm E(t,\bm x) |^2 \, \mathrm d \bm x$ exhibits the analytical decay rate approximately up to $t=35$. For rank $r=10$ the electrical energy starts to  deviate at the end of the time interval, while the solution with rank $r=15$ shows the correct rate within the full simulation time. 

Furthermore, it is known that the particle number, the total energy and the entropy, that is, the quantities
\begin{equation} \label{eq:invariants}
\begin{aligned}
&\int_{\Omega} f(t,\bm x,\bm v) \, \mathrm d \bm x \, \mathrm d \bm v, \quad
	\frac 1 2 \int_{\Omega} \vert \bm v\vert^2 \, f(t,\bm x,\bm v) \, \mathrm d \bm x \, \mathrm d \bm v 
		+ \frac 1 2 \int_{\Omega^{(\bm x)}} \vert \bm E(t,\bm x)\vert^2 \, \mathrm d \bm x, \\
&\int_{\Omega} \vert f(t,\bm x,\bm v)\vert^2 \, \mathrm d \bm x \, \mathrm d \bm v,
\end{aligned}
\end{equation}
are invariants of the exact solution. In Figure~\ref{fig:landau_fixed} we see that the mass and the total energy are almost conserved, whereas the entropy is only preserved up to a small error. 

\begin{figure}[h!]
\begin{center}
\includegraphics[width=0.95\textwidth]{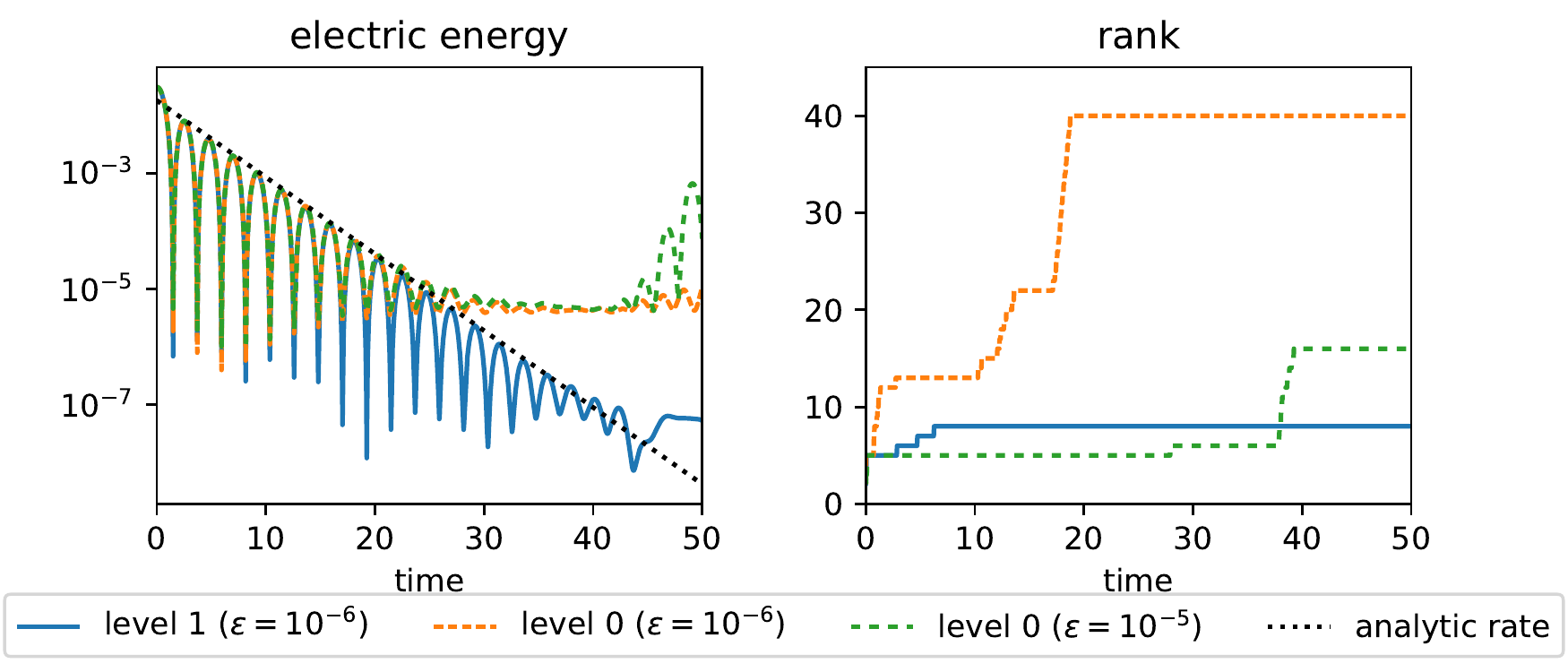}
\end{center}
\caption{Simulation results of the 2+2-dimensional Landau damping using the rank adaptive algorithm (Algorithm~\ref{alg:rauc}) for different tolerances $\epsilon$ and discretizations; electric energy including the analytical decay rate (left) and ranks (right). }
\label{fig:landau_rauc}
\end{figure}

For the rank adaptive case (Algorithm~\ref{alg:rauc}) the system was simulated for the same discretization (level 0) and different tolerances $\epsilon$. The corresponding results in Figure~\ref{fig:landau_rauc} show that the electric energy deteriorates at around $t=25$ for both tolerances $\epsilon=10^{-5}, 10^{-6}$. Although the rank increases up to $40$ (the maximal rank allowed) for the case of $\epsilon=10^{-6}$ the accuracy in the electric energy does not improve. To improve accuracy the simulation is carried out on a uniformly refined spatial and velocity mesh (level $1$). The results for a time step of $\Delta t = 0.00125$ and a tolerance of $\epsilon=10^{-6}$ are shown in Fig.~\ref{fig:landau_rauc}. The electric energy shows the correct decay in electric energy up to approximately $t=40$.

\begin{figure}[t]
\begin{center}
\includegraphics[width=0.9\textwidth,trim={0 4cm 0 14cm},clip]{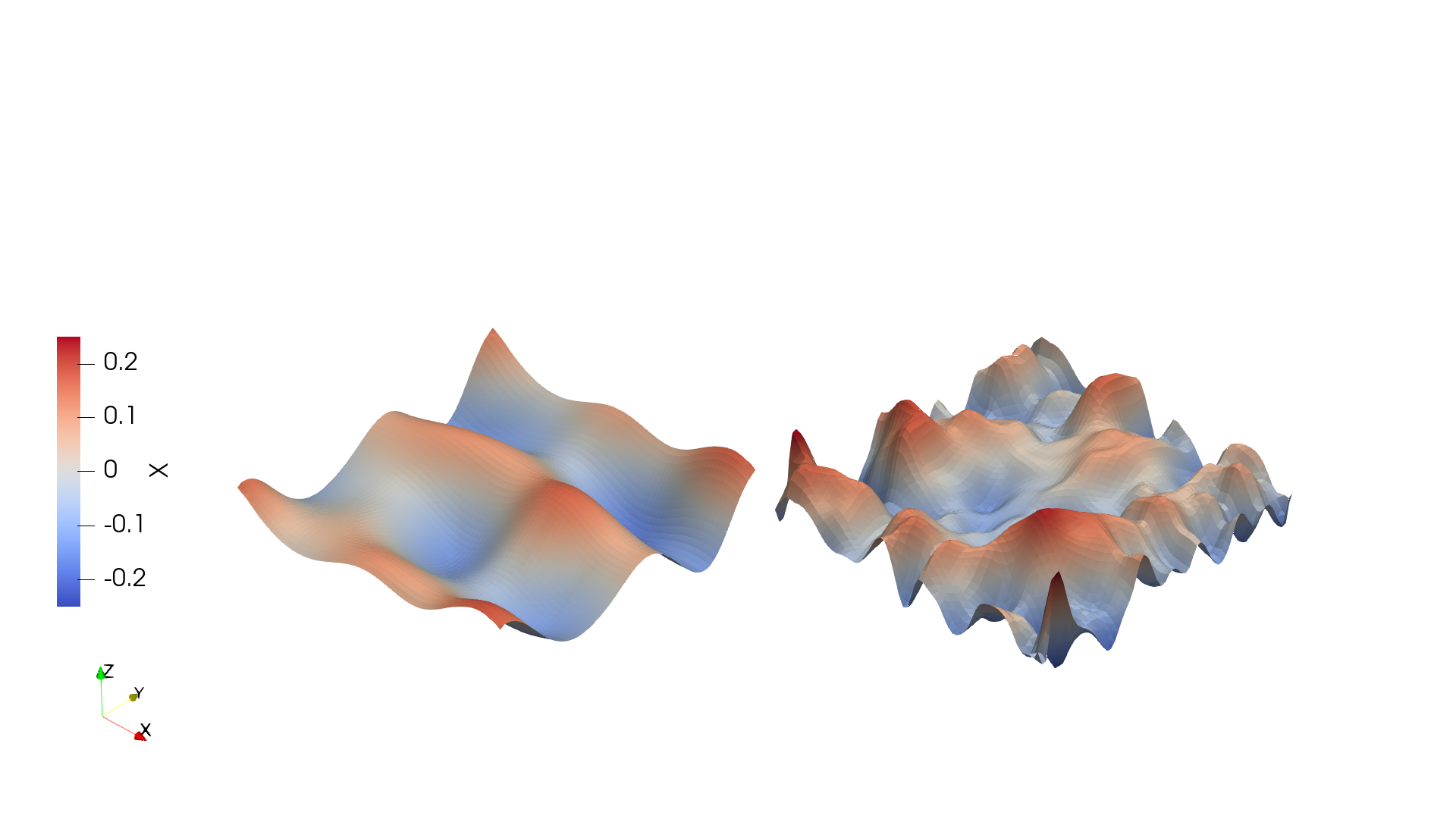}
\end{center}
\caption{Spatial basis functions $X_i(\bm x)$ at time $t=50$ for fixed rank ($r=15$, left) and rank adaptive (level 0, $\epsilon=10^{-5}$, right) simulation. }
\label{fig:X_landau}
\end{figure}

Investigating the computed bases $X$ and $V$ in more detail shows that in the rank adaptive case spurious oscillations are present. In Figure~\ref{fig:X_landau} two exemplary basis functions $X_i(\bm x)$ of level $0$ at time $t=50$  are displayed. In the fixed rank case ($r=15$) the basis function is much smoother than in the rank adaptive simulation ($\epsilon=10^{-5}$).

In summary spurious modes may enter in the course of the simulation for the rank adaptive simulation. However, the accuracy can be improved by refining the discretization. In contrast, the algorithm using a fixed rank seems to have a regularizing effect. It remains to future work to investigate this effect more closely.

\subsection{Inflow boundary condition with constant electrical field}

In the second example we focus on the boundary condition and compare the numerical to an analytical solution. In this setting, we are solving the transport equation~\eqref{eq:transport} in $2+2$ dimensions with a constant electrical field $\bm E = [0, 4]^T$ on a triangular spatial domain \[
\Omega_x = \{ (x_1,x_2) \,|\, -0.5 < x_1 < 0.5, \, -x_1/2+1/4 < x_2 < x_1/2 - 1/4 \}
\] 
with initial condition $f(0,\bm x, \bm v) = 0$ and $\Omega_v = \mathbb R^2$. The inflow 
\begin{equation} \label{eq:inflow_boltzmann}
f(t,\bm x, \bm v) = \bar f(t, \bm x, \bm v) \quad \textnormal{on } \Gamma^-. 
\end{equation}
on the boundary of $\Omega_x$ will be determined by a function $\bar f$, which is a solution of the same equation as for $f$, but on the whole domain $\mathbb R^2$ and with initial condition
\[ \bar f(0, \bm x, \bm v) = \bar f_0(\bm x, \bm v). \]
Here, $\bar f_0$ has compact spatial support just outside the triangular domain $\Omega_x$ and a compactly supported velocity distribution centred around $\bm v = [2,0]^T$. More precisely, we set
\begin{equation}\label{eq:definitionf_0}
\bar f_0(\bm x, \bm v) = 
    \phi\Big(\frac{x_1 - 0.5 - \sigma_x}{\sigma_x}\Big) \cdot 
    \phi\Big(\frac{x_2 - 0.1}{\sigma_x}\Big) \cdot 
    \phi\Big(\frac{v_1 - 2}{\sigma_v}\Big) \cdot 
    \phi\Big(\frac{v_2}{\sigma_v}\Big), 
\end{equation}
where $\sigma_x = 0.2$, $\sigma_v = 0.5$, and
\[ \phi(z) = \begin{cases} z^2 \cdot (2 |z| - 3) + 1 & |z| \le 1 \\ 0 & |z|>1 \end{cases}\]
is a $C^1$ function supported in $[-1,1]$ and centered around $0$. In consequence $\bar f_0$ is a product function supported in a four-dimensional cube with side lengths controlled by $\sigma_x$ and $\sigma_v$.

By the method of characteristics we obtain 
\begin{equation} \label{eq:analytical_solution}
\bar f(t, \bm x, \bm v) = \bar f_0(\bm x - \bm v \cdot t - \bm E \cdot t^2/2, \bm v + \bm E \cdot t). 
\end{equation}
Restricting $\bar f$ on $\Omega_x \times \mathbb R^2$ solves the original problem and can be used to assess the quality of its numerical solution $f^h_r$.

In order to use $\bar f$ as the inflow function in our dynamical low-rank approach, we have to work with an approximation in separable form instead. For that purpose, we use the fact that, by~\eqref{eq:analytical_solution} and our choice~\eqref{eq:definitionf_0}, $\bar f$ can be written as
\[ 
\bar f(t, \bm x, \bm v) = g_1(t,x_1, v_1) \cdot g_2(t,x_2, v_2). 
\]
The factors $g_1$ and $g_2$ are then evaluated for a given time on a regular fine grid $(x_i, v_i)$ and  singular value decompositions are computed. The product of the truncated SVDs is then used for computing the inflow~\eqref{eq:inflow_boltzmann} as well as the error of the numerical solution. In our experiments we used the truncation ranks of 25.

The transport equation is solved numerically on the time domain $[0,0.5]$. On the coarsest scale (level 0) we use a conforming triangulation of $\Omega_x$ with $339$ vertices. The velocity domain is chosen as $\Omega_v = [-4,4]^2$ with periodic boundary conditions and a regular triangulation with $4096$ vertices. For the numerical solution the rank adaptive algorithm (Algorithm~\ref{alg:rauc}) is used with a time step of $\Delta t = 0.005$, a truncation threshold of $\epsilon = 10^{-3}$, and a stabilization parameter $\delta = 10^{-2}$. Higher levels $\ell=1,2$ are obtained by uniformly refining the meshes in $\Omega_x$ as well as~$\Omega_v$. For these levels the parameters $\Delta t$ and $\epsilon$ are scaled by $2^{-2\ell}$.

\begin{figure}[t]
\begin{center}
$t=0.075$ \\
\includegraphics[width=0.38\textwidth]{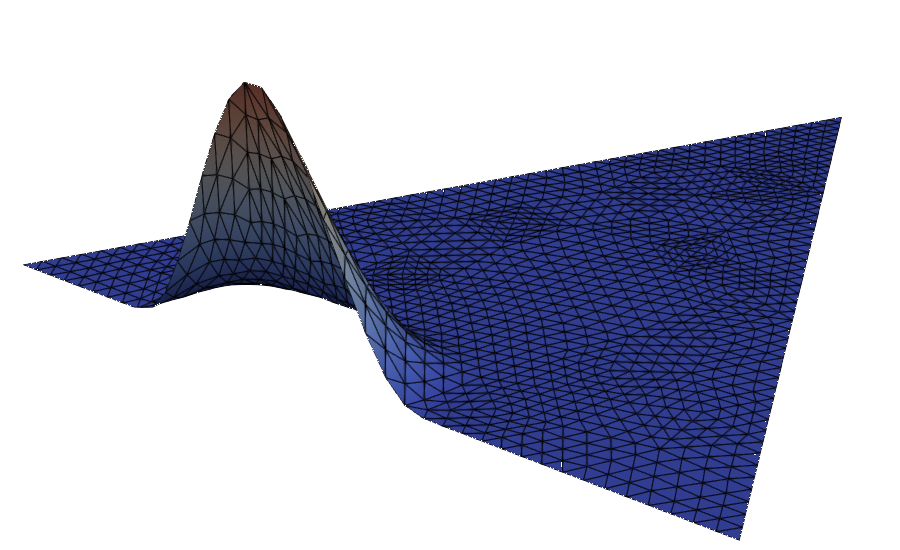} \hfill
\includegraphics[width=0.38\textwidth]{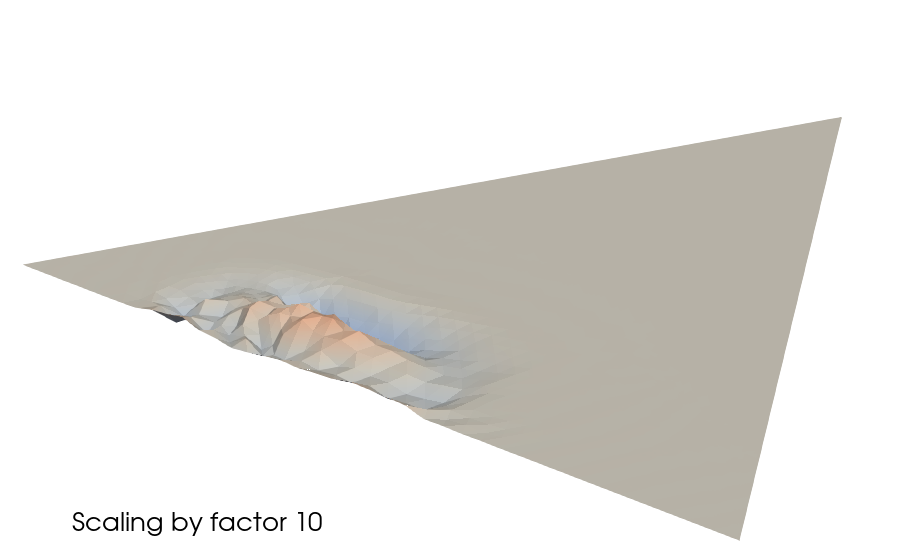} \\
$t=0.25$ \\
\includegraphics[width=0.38\textwidth]{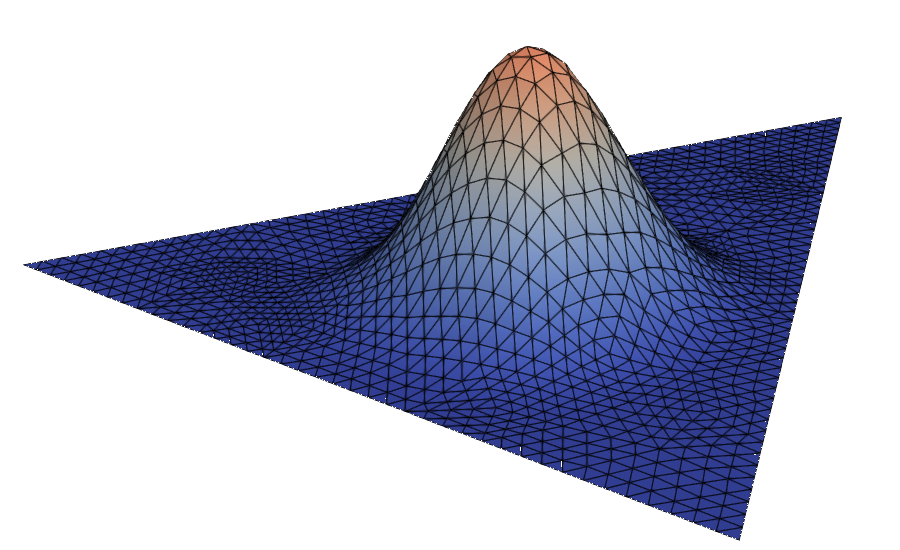} \hfill
\includegraphics[width=0.38\textwidth]{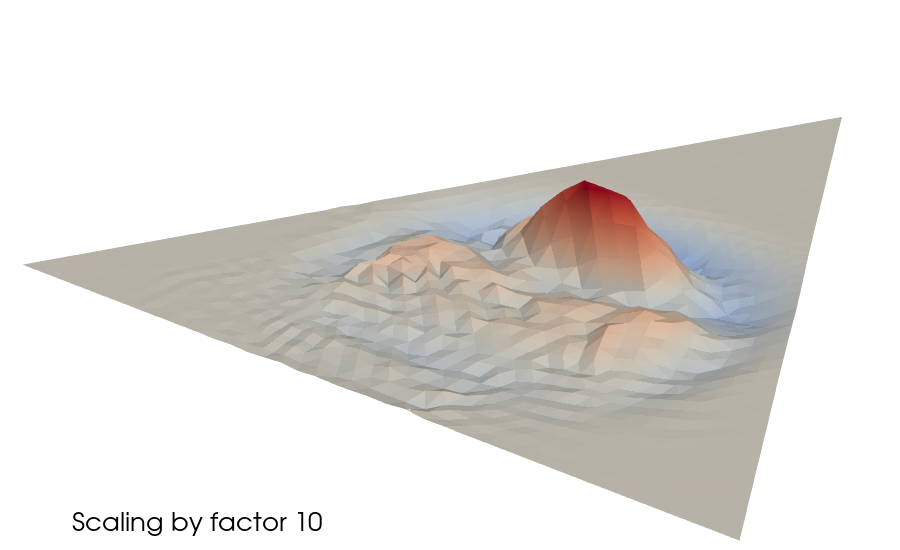} \\
$t=0.375$ \\
\includegraphics[width=0.38\textwidth]{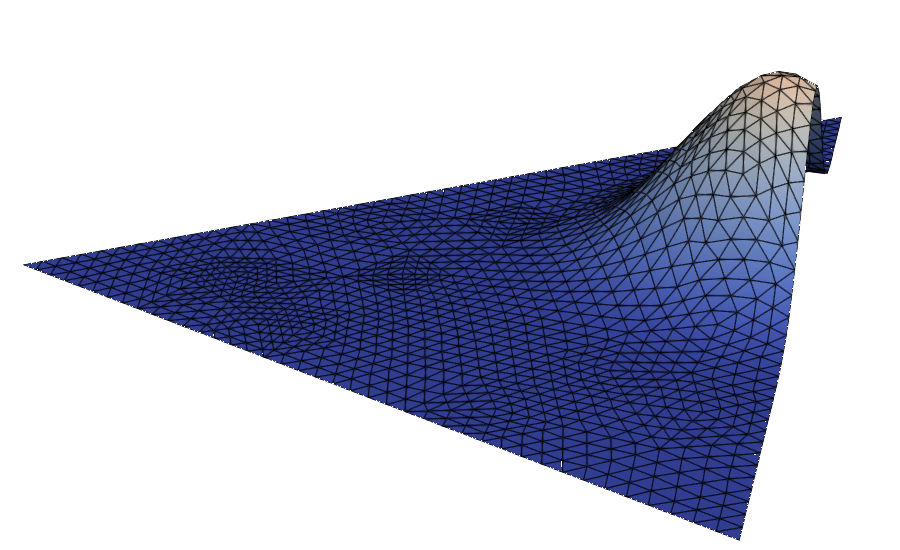} \hfill
\includegraphics[width=0.38\textwidth]{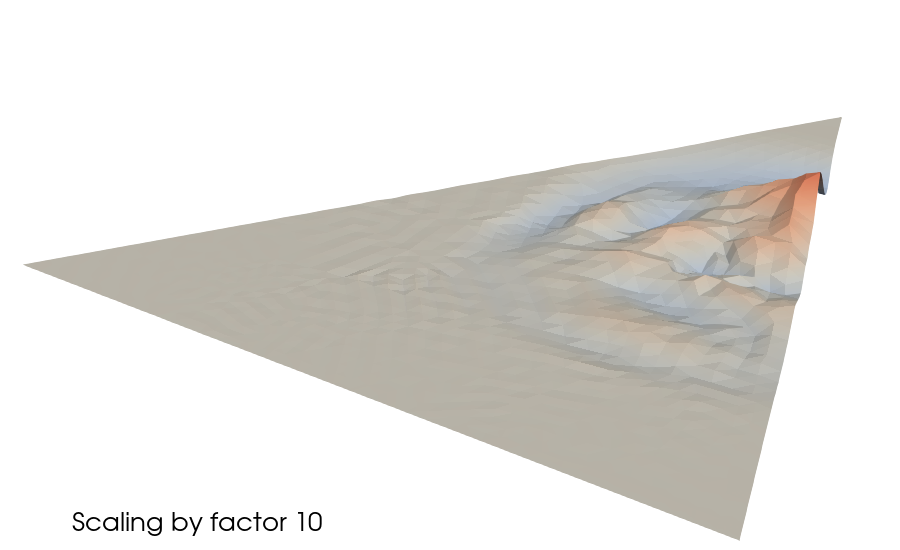} \\
\includegraphics[width=0.38\textwidth]{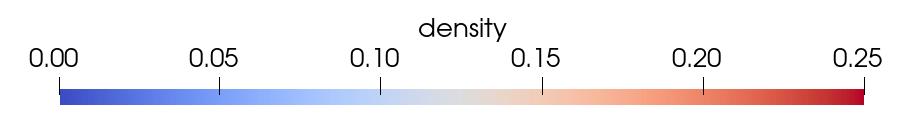} \hfill
\includegraphics[width=0.38\textwidth]{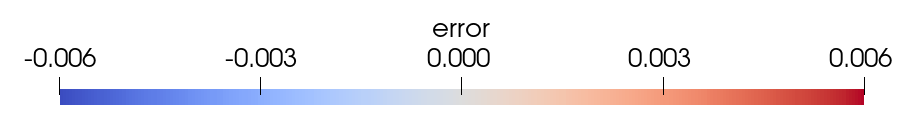}
\end{center}
\caption{Spatial density $\rho^h(t, \bm x)$ (left), see \eqref{eq:rhox}, and error scaled by a factor of 10 (right) of the numerical solution of \eqref{eq:transport} with constant electrical field for level $\ell=1$ at times $t = 0.075, \, 0.25, \, 0.375$.}
\label{fig:boltzmann_rhox}
\end{figure}

Figure~\ref{fig:boltzmann_rhox} shows the evolution of the computed spatial density 
\begin{equation} \label{eq:rhox}
\rho^h(t, \bm x) = \int_{\Omega_v} f_{r}^h(t,\bm x, \bm v) \, \mathrm d \bm v
\end{equation}
for level $\ell=1$ at different time steps together with its error computed using the analytical solution~$\bar f$.  At the beginning of the simulation, the density flows into the domain, is transported and finally leaves the domain. The error remains reasonably small during the whole process.

\begin{figure}[t]
\begin{center}
\vspace*{-4ex}
\includegraphics[width=0.75\textwidth]{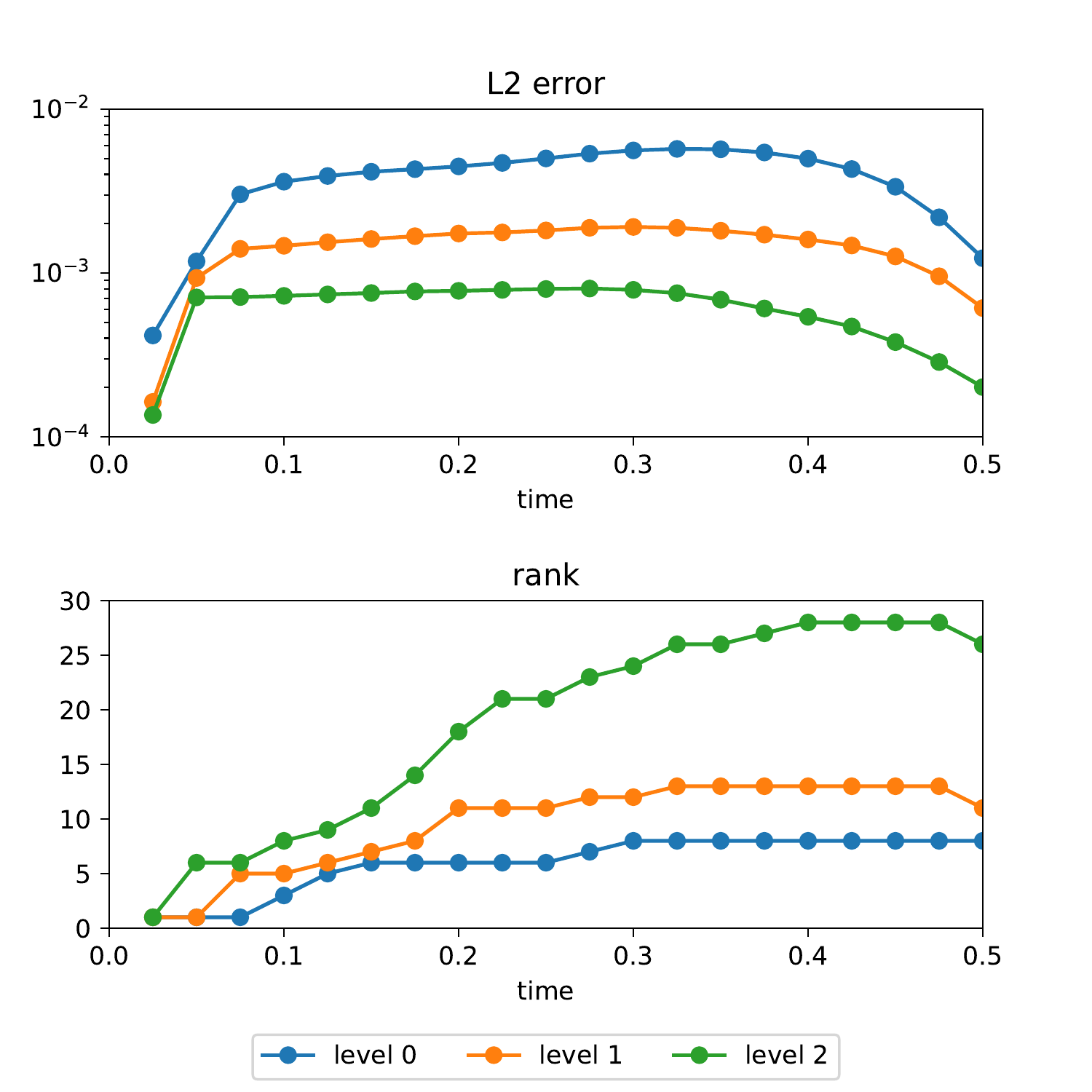}
\end{center}
\caption{Numerical solution of \eqref{eq:transport} for different levels of discretization. The upper graph shows the $L_2$ error computed with respect to the analytical solution $\bar f$, see \eqref{eq:analytical_solution},  on a uniformly refined grid. The lower graph shows the ranks used by the rank adaptive integrator. }
\label{fig:boltzmann}
\end{figure}

A more detailed investigation of the $L_2$ error is depicted in the upper part of Figure~\ref{fig:boltzmann}. It is obtained by comparing the numerical solutions $f_{r}^h$ to the analytical solution $\bar f$ on a once uniformly refined grid. As one can see, in the beginning the error increases for all levels but decays again for $t \ge 0.35$. At about that time a significant fraction of the density has left the domain already (see Figure~\ref{fig:boltzmann_rhox}). The maximal error decays as the level increases.

The lower part of the figure shows the ranks which were used by the rank adaptive algorithm. As more particles enter the domain and the distribution spreads out, the ranks increase. For higher levels of discretization a smaller truncation parameter is used in order to balance the low-rank approximation error with the discretization error, leading to higher ranks.

\section{Conclusion and outlook}

In this paper we have studied how to incorporate inflow boundary conditions in the dynamical low-rank approximation (DLRA) for the Vlasov--Poisson equation based on its weak formulation. The single steps in the projector splitting integrator, or the rank-adaptive unconventional integrator, can be interpreted as restrictions of the weak formulation to certain subspaces of the tangent space. The efficient solution of these sub-steps requires the separability of the boundary integrals, which is ensured for piecewise linear boundaries together with a separable inflow function. The resulting equations can be solved using FEM solvers based on Galerkin discretization. We confirmed the feasibility of our approach in numerical experiments.

As a next step, the conservation of physical invariants such as mass and momentum in the numerical schemes could be addressed, perhaps by extending methods from~\cite{Einkemmer2019,EINKEMMER2021110495,Einkemmer2022a,GuoWei2022,guo2022}. Enforcing nonnegativity of the density function in the DLRA approach is another open issue.

A potential advantage of the weak formulation for the sub-problems in the projector splitting integrator is that in principle it should allow for a great flexibility regarding the discretization spaces. In particular, they do not need to be fixed in advance and mesh-adaptive methods could be used for solving the sub-steps, provided suitable interpolation and prolongation operations are available, as well as adaptive solvers for Friedrichs' systems. We leave this as possible future work. Other elaborate adaptive integrators have been proposed in~\cite{hauck2022}, which also could be combined with our approach.

\paragraph*{Acknowledgements}
The work of A.U.~was supported by the Deutsche Forschungsgemeinschaft (DFG, German Research Foundation) – Projektnummer 506561557.

\paragraph*{Version of record}
This version of the article has been accepted for publication, after peer review but is not the version of record and does not reflect post-acceptance improvements, or any corrections. The version of record is available online at: \url{https://doi.org/10.1007/s10543-024-01019-8}.

\end{document}